\documentclass[11pt]{amsart}
\usepackage[english]{babel}
\usepackage{egothic}
\usepackage[T1]{fontenc}

\usepackage{amsmath}
\usepackage{amsthm}
\usepackage{amssymb}
\usepackage[all]{xy}
\usepackage{mathrsfs}
\usepackage{enumerate}
\usepackage{physics}
\usepackage[notcite, final, notref]{showkeys}
\usepackage{dsfont}
\usepackage[dvips]{color}
\newtheorem{theorem}{Theorem}[section]

\newtheorem{lemma}[theorem]{Lemma}
\newtheorem{proposition}[theorem]{Proposition}
\newtheorem{corollary}[theorem]{Corollary}
\newtheorem{definition}[theorem]{Definition}

\theoremstyle{definition}
\newtheorem{remark}[theorem]{Remark}

\newtheorem{example}[theorem]{Example}

\newcommand{\w}{\omega}

\renewcommand{\Re}{\mathrm{Re}}

\newcommand{\z}{\overline{z}}

\usepackage{xcolor}

\title[]{Reproducing kernel Hilbert spaces of polyanalytic functions of infinite order}

\author[D. Alpay]{Daniel Alpay}
\address{(DA) Schmid College of Science and Technology \\
Chapman University\\
One University Drive
Orange, California 92866\\
USA}
\email{alpay@chapman.edu}

\author[F. Colombo]{Fabrizio Colombo}
\address{(FC) Politecnico di
Milano\\Dipartimento di Matematica\\Via E. Bonardi, 9\\20133 Milano\\Italy}
\email{fabrizio.colombo@polimi.it}

\author[K. Diki]{Kamal Diki}
\address{(KD) Schmid College of Science and Technology \\
Chapman University\\
One University Drive
Orange, California 92866\\
USA}
\email{diki@chapman.edu}

\author[I. Sabadini]{Irene Sabadini}
\address{(IS) Politecnico di
Milano\\Dipartimento di Matematica\\Via E. Bonardi, 9\\20133 Milano\\Italy}
\email{irene.sabadini@polimi.it}

\begin{document}
\maketitle
\tableofcontents
\begin{abstract}
In this paper we introduce reproducing kernel Hilbert spaces of polyanalytic functions of infinite order. First we study in details the counterpart of the Fock space and related results in this framework. In this case the kernel function is given by $\displaystyle e^{z\overline{w}+\overline{z}w}$ which can be connected to kernels of polyanalytic Fock spaces of finite order. Segal-Bargmann and Berezin type transforms are also considered in this setting. Then, we study the reproducing kernel Hilbert spaces of complex-valued functions with reproducing kernel $\displaystyle\frac{1}{(1-z\overline{w})(1-\overline{z}w)}$ and $\displaystyle\frac{1}{1-2{\rm Re}\, z\overline{w}}$. The corresponding backward shift operators are introduced and investigated.
\end{abstract}


\section{Introduction and preliminary results}
\setcounter{equation}{0}
Polyanalytic functions were introduced in 1908 by Kolossov to solve problems in elasticity theory, see \cite{kolossov}. For a general introduction to this topic see \cite{Alpay2015, Balk1991, balk_ency}. In more recent times, this function theory was
studied by several authors from different perspectives, see
\cite{abreu, abreufeicht, agranovsky, begehr, vasilevski} and the references
therein.
Polyanalytic functions are used also to study sampling and interpolation problems on Fock spaces using time frequency analysis techniques such as short-time Fourier transform (STFT) or Gabor transforms, see  \cite{A2010}. In the next parts of this introduction, we collect some basic definitions and explain what we mean by polyanalytic functions of infinite order in our setting. Some important facts that will be needed in the sequel will be also revised. Then, we will explain the general construction of the kernels associated to the reproducing kernel Hilbert spaces of polyanalytic functions of infinite order, which will be studied in this paper. We conclude by describing the contents of the paper.
\subsection{Definitions}
A complex valued function $f:\Omega\subset \mathbb{C}\longrightarrow \mathbb{C}$ which belongs to the kernel of a power $n\geq 1$ of the classical Cauchy-Riemann operator $\displaystyle \frac{\partial}{\partial \overline{z}}$, that is $$\displaystyle \frac{\partial^n}{\partial \overline{z}^n}f(z)=0, \quad \forall  z\in\Omega,$$
 is called a polyanalytic function of order $n$.
\smallskip
An interesting fact regarding these functions is that any polyanalytic function of order $n$ can be decomposed in terms of $n$ analytic functions so that we have a decomposition of the following form
\begin{equation}
f(z)=\displaystyle \sum_{k=0}^{n-1}\overline{z}^kf_k(z),
\end{equation}
for which all $f_k$ are analytic functions on $\Omega$. In particular, expanding each analytic component using the series expansion theorem lead to an expression of this form

\begin{equation}\label{exp1}
f(z)=\displaystyle \sum_{k=0}^{n-1}\sum_{j=0}^{\infty}\overline{z}^kz^ja_{k,j},
\end{equation}
where $(a_{k,j})$ are complex coefficients.
In this paper, we are interested by the case where the expansion \eqref{exp1} is of infinite order, which means that we consider functions of the form
\begin{equation}\label{exp2}
f(z)=\displaystyle \sum_{k=0}^{\infty}\sum_{j=0}^{\infty}\overline{z}^kz^ja_{k,j},
\end{equation}

which will be called polyanalytic functions of infinite order. We note that such functions were discussed in \cite{Balk1991, balk_ency} in which they were mentioned as conjugate analytic functions.
\\ \\
For $n=1,2,...$ we recall that polyanalytic Fock spaces of order $n$ can be defined as follows
$$\mathcal{F}_n(\mathbb{C}):=\left\lbrace g\in H_n(\mathbb{C}), \quad \frac{1}{\pi}\int_{\mathbb{C}}|g(z)|^2e^{-|z|^2}dA(z)<\infty \right\rbrace.$$
The reproducing kernel associated to the space $\mathcal{F}_n(\mathbb{C})$ is given by
\begin{equation}\label{Kn}
K_n(z,w)=e^{z\overline{w}}\displaystyle \sum_{k=0}^{n-1}\frac{(-1)^k}{k!}{n \choose k+1}|z-w|^{2k},
\end{equation}
for every $z,w\in\mathbb{C}.$
\subsection{The kernels construction: general discussion}
Consider a function $F(z_1,z_2)$ in the Hardy space of the bidisk $\mathbb{D}^2$ in $\mathbb{C}^2$. Then, $g(z)=F(z,z)$ belongs to the Bergman space of the disk $\mathbb{D}$, and the map $F\mapsto g$ is onto and contractive, but not one-to-one.
For instance, the polynomials $z_1^nz_2^m$ belong to $\mathbf H^2(\mathbb D^2)$ and have the same image $z^{s}$ with $n+m=s$. On the other hand the map $g(z)=F(z,\z)$ is one-to-one, and its image is the
reproducing kernel  Hilbert space with reproducing kernel
\[
  \frac{1}{(1-z\overline{w})(1-\z w)}
\]
The corresponding reproducing kernel Hilbert space consists of polyanalytic functions of infinite order. \\ \\ Motivated by the above discussion
we consider a function $c(\mathbf z, \mathbf w)$ positive definite in some open subset $\Omega$ of $\mathbb C^{2N}$, and analytic in $\mathbf z$ and $\overline{\mathbf w}$.
We assume that
\begin{equation}
  \Omega_s=\left\{z\in\mathbb C^N\,:\, (z,\z)\in\Omega\right\}
\label{sym}
\end{equation}
is open and non-empty. The function
\[
k(z,w)=c((z,\z),(\w,\overline{\w}))
\]
is then positive definite in $\Omega_s$. The purpose of the present work is to study the corresponding reproducing kernel Hilbert spaces of polyanalytic functions of infinite
order. For instance in the case of the Fock space with reproducing kernel $c(\z,\w)=e^{\sum_{n=1}^{2N}z_n\overline{\w_n}}$ we have the kernel
  \begin{equation}
k(z,w)=e^{\sum_{n=1}^N(z_n\overline{w_n}+\overline{z_n}w_n)}
\end{equation}
while in the case of the Drury-Arveson space with reproducing kernel $\frac{1}{1-\sum_{n=1}^{2N}z_n\overline{w_n}}$ which is positive definite in the open unit ball of
$\mathbb C^{2N}$, the corresponding kernel is
\begin{equation}
  \frac{1}{1-\sum_{n=1}^N(z_n\overline{w_n}+\overline{z_n}w_n)}
\end{equation}
positive definite in the open ball of $\mathbb C^N$ centered at the origin and with radius $\frac{1}{\sqrt{2}}$. In this paper we will focus on $N=1$.\\

A general family of examples correspond to
\[
k(z,w)=K_1(z,w)\overline{K_2(z,w)}
\]
where $K_1$ and $K_2$ are analytic kernels, or, in the matrix-valued case,
\[
k(z,w)=K_1(z,w)\otimes \overline{K_2(z,w)}.
\]

The structure of the paper is as follows: in Section 2 we introduce the kernel function $K$ associated to the polyanalytic Fock space $\mathcal{SF}(\mathbb C)$ of infinite order and we study various properties. We give a sequential characterization of the space $\mathcal{SF}(\mathbb C)$ and, in particular, we prove that the creation and annihilation operators are adjoint of each other. We also introduce and study two backward shift operators. In Section 3 we prove that by taking the power series of the polyanalytic Fock kernels of finite order $(K_n)_{n\geq 0}$ we obtain the kernel function $K$ multiplied up to an exponential kernel. In Section 4 and 5 we study Segal-Bargmann and Berezin type transforms and some related operators. In Section 6 we present the polyanalytic Hardy space of infinite order and we study the Gleason problem. We also prove some results on the backward shift operator in this setting.  Finally, Section 7 is devoted to the case of Drury-Arveson space.


\section{The polyanalytic Fock space of infinite order and associated kernel}
\setcounter{equation}{0}
We denote by $M_z$ and $M_{\overline{z}}$ the multiplication operators by $z$ and $\overline{z}$. Then, we will prove the following main result
\begin{theorem}
 The reproducing kernel Hilbert space with reproducing kernel $e^{z\overline{w}+\z w}$ is, up to a multiplicative positive factor, the only reproducing kernel
Hilbert space of polyanalytic functions of infinite order, regular at the origin, and for which
\begin{eqnarray}
\left(\frac{\partial}{\partial z}\right)^*&=&M_z\\
\left(\frac{\partial}{\partial \overline{z}}\right)^*&=&M_{\overline{z}}.
 \end{eqnarray}
\end{theorem}
To this end, we need the following:
\begin{definition}
We consider the kernel function given by
\begin{equation}\label{newFock-kernel}
K(z,w)=e^{z\overline{w}+\z w}=e^{2Re(z\overline{w})}, \quad \forall (z,w)\in\mathbb{C}^2.
\end{equation}
We denote by $(\mathcal{H}(K),\langle \cdot, \cdot \rangle_{\mathcal{H}(K)})$ the reproducing kernel
Hilbert space associated to the kernel function \eqref{newFock-kernel}.
\end{definition}

\begin{proposition}
The function $K:\mathbb{C}\times \mathbb{C}\longrightarrow \mathbb{C}$ defined by \eqref{newFock-kernel} is a positive definite kernel.
\end{proposition}
\begin{proof}
 It is clear that we have
\begin{equation}\label{FFbar}
K(z,w)=F(z,w)\overline{F(z,w)},
\end{equation}

for every $z,w\in\mathbb{C}$ and $F$ denotes the reproducing kernel of the classical Fock space $\mathcal{F}(\mathbb{C})$. Thus, since $\overline{F(z,w)}$ is also a positive definite kernel we can conclude, since $K$ is a product of positive definite kernels.
\end{proof}

We observe that the following integral representation holds
\begin{proposition}
It holds that \begin{equation}
\displaystyle \frac{1}{\pi}\int_{\mathbb{C}}K(z,w)e^{-|w|^2}dA(w)=e^{|z|^2}=\sqrt{K(z,z)}, \text{ for any } z\in\mathbb{C}.
\end{equation}

\end{proposition}
\begin{proof}
We set $w=x+iy$, we identify $\mathbb{C}$ with  $\mathbb{R}^2$ and use the classical Gaussian integral
$$\displaystyle\int_{\mathbb{R}} e^{-at^2+bt}dt=\sqrt{\frac{\pi}{a}}e^{\frac{b^2}{4a}}, \quad a>0, b\in\mathbb{C}.$$
We have

\[
    \begin{split}
   \displaystyle   \frac{1}{\pi}\int_{\mathbb{C}}K(z,w)e^{-|w|^2}dA(w) &= \frac{1}{\pi}\int_{\mathbb{C}}e^{z\overline{w}+\z w}e^{-|w|^2}dA(w) \\
      &=  \frac{1}{\pi}\int_{\mathbb{C}}e^{z(x-iy)+\z(x+iy)}e^{-(x^2+y^2)}dxdy \\
        &= \frac{1}{\pi}\left(\int_{\mathbb{R}}e^{x(\overline{z}+z)-x^2}dx\right)\left(\int_{\mathbb{R}}e^{y i(\overline{z}-z)-y^2}dy \right)
        \\
        &= e^{\frac{(z+\overline{z})^2}{4}}e^{-\frac{(\overline{z}-z)^2}{4}}
        \\
        &= e^{|z|^2}
         \\
        &=\sqrt{K(z,z)},
          \end{split}
   \]
   as stated.
\end{proof}
Next result is simple but very useful:
\begin{proposition} \label{Fockpptd}
For any $z,w\in \mathbb{C}$ it holds that
\begin{itemize}
\item[i)] $\displaystyle \frac{\partial}{\partial z}K(z,w)=\overline{w}K(z,w)$ and $\displaystyle \frac{\partial}{\partial \overline{z}}K(z,w)=wK(z,w)$.
\item[ii)] $\displaystyle \frac{\partial}{\partial w}K(z,w)=\overline{z}K(z,w)$ and $\displaystyle\frac{\partial}{\partial \overline{w}}K(z,w)=zK(z,w)$.
\end{itemize}
\end{proposition}
\begin{proof}
It is an easy calculation based on the expression of the kernel function $K(z,w)$ given by formula \eqref{newFock-kernel}.
\end{proof}
\begin{corollary}\label{Fockppt}
For any $z,w\in \mathbb{C}$ and $n=1,2,...$ it holds that
\begin{itemize}
\item[i)] $\displaystyle \frac{\partial^n}{\partial z^n}K(z,w)=\overline{w}^nK(z,w)$ and $\displaystyle \frac{\partial^n}{\partial \overline{z}^n}K(z,w)=w^nK(z,w)$.
\item[ii] $\displaystyle \frac{\partial^n}{\partial w^n}K(z,w)=\overline{z}^nK(z,w)$ and $\displaystyle\frac{\partial^n}{\partial \overline{w}^n}K(z,w)=z^nK(z,w)$.
\end{itemize}
\end{corollary}
\begin{proof}
We will prove i), the other statements follow similar arguments. We apply Proposition \ref{Fockpptd} and get $$\displaystyle \frac{\partial}{\partial z}K(z,w)=\overline{w}K(z,w).$$
Thus, we apply a second time the complex derivative with respect to the variable $z$, we use again Proposition \ref{Fockpptd} and get $$\displaystyle \frac{\partial^2}{\partial z^2}K(z,w)=\overline{w}^2K(z,w).$$
Then, we repeat the same calculation $n$-times and obtain $$\displaystyle \frac{\partial^n}{\partial z^n}K(z,w)=\overline{w}^nK(z,w).$$
\end{proof}
Now, let us consider the commutator operators given by \begin{equation}
\displaystyle \left[\frac{\partial}{\partial z}, M_z\right]:=\frac{\partial}{\partial z}M_z-M_z\frac{\partial}{\partial z}
\end{equation}
and

\begin{equation}
\displaystyle \left[\frac{\partial}{\partial \overline{z}}, M_{\overline{z}}\right]:=\frac{\partial}{\partial  \overline{z}}M_{\overline{z}}-M_{\overline{z}}\frac{\partial}{\partial \overline{z}}.
\end{equation}
Then, we can prove the following
\begin{proposition}\label{Commutator}
For any $z,w\in\mathbb{C}$ we have
$$\displaystyle \left[\frac{\partial}{\partial z}, M_z\right]K(z,w)=K(z,w)$$
and
$$\displaystyle \left[\frac{\partial}{\partial \overline{z}}, M_{\overline{z}}\right]K(z,w)=K(z,w).$$

\end{proposition}
\begin{proof}
We have
  \[
    \begin{split}
      \displaystyle\frac{\partial}{\partial z}M_z K(z,w)&=\frac{\partial}{\partial z}\left(zK(z,w) \right)\\
      &=z\frac{\partial}{\partial z}K(z,w)+K(z,w)\\
        &=z\overline{w}K(z,w)+K(z,w).
          \end{split}
   \]
On the other hand
 \[
    \begin{split}
      \displaystyle M_z\frac{\partial}{\partial z}K(z,w)&=M_z\left(\overline{w}K(z,w)\right)\\
      &=z\overline{w}K(z,w),\\
          \end{split}
   \]
   hence, we obtain $$\displaystyle \left[\frac{\partial}{\partial z}, M_z\right]K(z,w)=K(z,w).$$

   In a similar way we can prove that
   $$\displaystyle \left[\frac{\partial}{\partial \overline{z}}, M_{\overline{z}}\right]K(z,w)=K(z,w).$$

\end{proof}
Thanks to the reproducing kernel property we have this more general result
\begin{theorem}
For any $f\in\mathcal{H}(K)$, the following identities hold
$$\displaystyle \left[\frac{\partial}{\partial z}, M_z\right]f=f$$
and
$$\displaystyle \left[\frac{\partial}{\partial \overline{z}}, M_{\overline{z}}\right]f=f.$$
\end{theorem}
\begin{proof}
Let $f\in\mathcal{H}(K)$, we know that $$f(z)=<f,K_z>_{\mathcal{H}(K)}, \text{ for any } z\in \Omega.$$
Thus, for any $z\in \Omega$ we apply Proposition \ref{Commutator} and get    \[
    \begin{split}
      \displaystyle  \left[\frac{\partial}{\partial z}, M_z\right]f(z) &= <f, \left[\frac{\partial}{\partial z}, M_z\right]K_z>_{\mathcal{H}(K)}\\
      &= <f,K_z>_{\mathcal{H}(K)}\\
        &=f(z).
          \end{split}
   \]
   Hence, it follows that $$\displaystyle \left[\frac{\partial}{\partial z}, M_z\right]f=f.$$
   In the same way we can prove that $$\displaystyle \left[\frac{\partial}{\partial \overline{z}}, M_{\overline{z}}\right]f=f.$$
\end{proof}
As a consequence of Proposition \ref{Fockpptd} we can study more properties of the kernel function in \eqref{newFock-kernel}.
\begin{proposition}\label{deltadelta}
For any $z, w\in\mathbb{C}$ we have
\begin{equation}
\displaystyle \Delta_z K(z,w)=4|w|^2K(z,w)
\end{equation}
and \begin{equation}
\Delta_w\Delta_zK(z,w)=16\left(1+w\overline{z}+\overline{w}z+|w|^2|z|^2\right)K(z,w)=16|1+\overline{w}z|^2K(z,w).
\end{equation}
\end{proposition}
\begin{proof}
We note that using  \eqref{newFock-kernel} and Proposition \ref{Fockpptd} we have
\[
    \begin{split}
      \displaystyle \frac{\partial^2}{\partial \overline{z}\partial z}K(z,w) &= \overline{w}\frac{\partial}{\partial \overline{z}}K(z,w) \\
      &= \overline{w}wK(z,w)\\
        &=|w|^2K(z,w).
          \end{split}
   \]
   Since $\displaystyle \frac{\partial^2}{\partial \overline{z}\partial z}=\frac{1}{4}\Delta_z$ we conclude that
  \begin{equation}\label{deltazK}
  \displaystyle \Delta_z K(z,w)=4|w|^2K(z,w).
  \end{equation}

Furthermore, by taking the derivative  with respect to $w$ we get \[
    \begin{split}
      \displaystyle \frac{\partial}{\partial w}\Delta_z K(z,w) &= 4\overline{w}\frac{\partial}{\partial w}\left(wK(z,w)\right) \\
      &= 4\overline{w}(K(z,w)+w\overline{z}K(z,w))\\
        &=4(1+w\overline{z})\overline{w}K(z,w).
          \end{split}
   \]
Then, we apply the derivative with respect to $\overline{w}$, develop the computations and get  \[
    \begin{split}
      \displaystyle \frac{\partial^2}{\partial \overline{w}\partial w}\Delta_z K(z,w) &= 4(1+w\overline{z})\frac{\partial}{\partial \overline{w}}\left(\overline{w} K(z,w)\right) \\
      &=4(1+w\overline{z})(1+\overline{w}z)K(z,w) \\
        &=4|1+\overline{w}z|^2K(z,w),
          \end{split}
   \]

   from which we conclude that $$\Delta_w\Delta_zK(z,w)=16|1+\overline{w}z|^2K(z,w).$$

\end{proof}
\begin{corollary}
Let $\Omega=\mathbb{D}$ denote the unit disk, we consider the operator $$T_{z,w}:=\frac{\Delta_w\Delta_z}{16|1+\overline{w}z|^2}, \quad \forall (z,w)\in\Omega\times \Omega.$$
Then, we have \begin{equation}
T_{z,w}K(z,w)=K(z,w), \quad \forall (z,w)\in\Omega\times \Omega.
\end{equation}
\end{corollary}
\begin{proof}
It is a direct consequence of Proposition \ref{deltadelta}.
\end{proof}
\begin{remark}
From \eqref{deltazK} we deduce that if $w$ is a fixed parameter, then the kernel function $K(z,w)$ can be seen as an eigenfunction of the Laplace operator $\Delta_z$ with eigenvalue given by $\displaystyle 4|w|^2$.
\end{remark}
\begin{definition}
The polyanalytic Fock space $\mathcal{SF}(\mathbb C)$ of infinite order is the set of functions of the form \begin{equation}
f(z)=\displaystyle \sum_{n=0}^{\infty}\overline{z}^nf_n(z),
\end{equation}
satisfying the conditions
\begin{enumerate}
\item[i)] $f_n\in \mathcal{F}(\mathbb{C})$ for any $n\geq 0$;
\item[ii)]$\displaystyle ||f||^{2}_{\mathcal{SF}(\mathbb{C})}= \sum_{n=0}^{\infty} n! ||f_n||^{2}_{\mathcal{F}(\mathbb{C})}<\infty.$
\end{enumerate}
Then, we consider the scalar product on $\mathcal{SF}(\mathbb{C})$ given by \begin{equation}\displaystyle
\langle f,g \rangle_{\mathcal{SF}(\mathbb{C})}:=\sum_{k=0}^{\infty}k!\langle f_k,g_k \rangle_{\mathcal{F}(\mathbb{C})},
\end{equation}
 for any $f=\displaystyle \sum_{k=0}^{\infty}\overline{z}^kf_k$ and $g=\displaystyle \sum_{k=0}^{\infty}\overline{z}^kg_k$ with $f_k,g_k\in\mathcal{F}(\mathbb{C})$ for every $k\geq 0$.
\end{definition}

\begin{proposition}\label{seqNF}
A function $f:\mathbb{C}\longrightarrow \mathbb{C}$ belongs to $\mathcal{SF}(\mathbb{C})$ if and only if $f$ is of the form $$\displaystyle f(z)=\sum_{(m,n)\in\mathbb{N}^2}z^m\overline{z}^n\alpha_{m,n},$$
with $(\alpha_{m,n})\subset\mathbb{C}$ and such that \begin{equation}\displaystyle
||f||^{2}_{\mathcal{SF}(\mathbb{C})}=\sum_{(m,n)\in\mathbb{N}^2}m!n!|\alpha_{m,n}|^2<\infty.
\end{equation}
Moreover, if for any $(m,n)\in\mathbb{N}^2$ we set $\phi_{m,n}(z,\overline{z})=\displaystyle \frac{z^m\bar{z}^n}{\sqrt{m!n!}}$ then, the family of functions $\lbrace \phi_{m,n} \rbrace_{m,n\geq 0}$ form an orthonormal basis of $\mathcal{SF}(\mathbb{C})$.
\end{proposition}
\begin{proof}
Let
$$\displaystyle f(z)=\sum_{(m,n)\in\mathbb{N}^2}\overline{z}^nz^m\alpha_{m,n},$$
with $(\alpha_{n,m})\subset\mathbb{C}$.

Setting $f_n(z)=\displaystyle \sum_{m=0}^\infty z^m\alpha_{m,n}$, it is clear that $\displaystyle f(z)=\sum_{n=0}^{\infty}\overline{z}^nf_n(z)$. Moreover, we have

$$\displaystyle ||f||^{2}_{\mathcal{SF}(\mathbb{C})}=\sum_{n=0}^{\infty}n!||f_n||^{2}_{\mathcal{F}(\mathbb{C})}=\sum_{(m,n)\in\mathbb{N}^2}m!n!|\alpha_{m,n}|^2.$$
Therefore, $f$ belongs to the space $\mathcal{SF}(\mathbb{C})$ if and only if

$$\displaystyle
||f||^{2}_{\mathcal{SF}(\mathbb{C})}=\sum_{(m,n)\in\mathbb{N}^2}m!n!|\alpha_{m,n}|^2<\infty.$$
On the other hand, easy computations lead to $$\langle f,\phi_{m,n}\rangle_{\mathcal{SF}(\mathbb{C})}=\sqrt{n!m!}\alpha_{m,n}, \quad \forall m,n \geq 0.$$

If $\langle f,\phi_{m,n}\rangle_{\mathcal{SF}(\mathbb{C})}=0$ for any $m,n\geq 0$, then we have $\alpha_{m,n}=0$ for any $m,n\geq 0$.
We note also that $$\langle \phi_{m,n},\phi_{m,n}\rangle_{\mathcal{SF}(\mathbb{C})}=1 \text{ and }  \langle \phi_{m,n},\phi_{p,q}\rangle_{\mathcal{SF}(\mathbb{C})}=0 \text{ whenever } (m,n)\neq (p,q).$$
In particular, this shows that $\lbrace \phi_{m,n} \rbrace_{m,n\geq 0}$ form an orthonormal basis of $\mathcal{SF}(\mathbb{C})$. This ends the proof.
\end{proof}
\begin{example}
We recall the complex Hermite polynomials introduced in \cite{Ito}
\begin{equation}
H_{m,n}(z,\overline{z}):=\displaystyle \sum_{k=0}^{\min{(m,n)}}(-1)^k k! {m\choose k}{n \choose k}z^{m-k}\overline{z}^{n-k}.
\end{equation}
It is easy to prove that $H_{m,n}$ belong to $\mathcal{SF}(\mathbb{C})$.
\end{example}
We now provide a sequential charachterization of the space $\mathcal{H}(K)$
\begin{theorem}
We have
$$\mathcal{H}(K)=\mathcal{SF}(\mathbb{C}).$$
Moreover, it holds that

\begin{equation}
K(z,w)=\displaystyle \sum_{m,n=0}^{\infty}\phi_{m,n}(z,\overline{z})\overline{\phi_{m,n}(w,\overline{w})}, \quad \text{ for any } z,w\in\mathbb{C}.
\end{equation}
\end{theorem}
\begin{proof}
Since $(\phi_{m,n})_{m,n\geq 0}$ is an orthonormal basis of the space $\mathcal{SF}(\mathbb{C})$, the associated reproducing kernel is given by the convergent series
$$\displaystyle \sum_{m,n=0}^{\infty}\phi_{m,n}(z,\bar{z})\overline{\phi_{m,n}(w,\bar{w})}<\infty, \text{ for any } z,w\in\mathbb{C}.$$
More precisely, for any $(z,w)\in\mathbb{C}^2$ we have the equalities
\[
    \begin{split}
   \displaystyle  \sum_{m,n=0}^{\infty}\phi_{m,n}(z,\bar{z})\overline{\phi_{m,n}(w,\bar{w})}
      &= \sum_{m,n=0}^{\infty}\frac{z^m\bar{z}^n \bar{w}^m w^n}{m!n!} \\
        &= \left(\sum_{m=0}^{\infty} \frac{z^m\bar{w}^m}{m!}\right) \left(\sum_{n=0}^{\infty} \frac{w^n\bar{z}^n}{n!}\right)
        \\
        &= e^{z\bar{w}}e^{w\bar{z}}
        \\
        &=e^{z\bar{w}+w\bar{z}}
         \\
        &=K(z,w).
          \end{split}
   \]
\end{proof}

\begin{theorem}\label{Adj}
It holds that
\begin{equation}\displaystyle
\langle \frac{\partial}{\partial z}f,g\rangle_{\mathcal{SF}(\mathbb{C})}=\langle f,M_{z}g\rangle_{\mathcal{SF}(\mathbb{C})},
\end{equation}
moreover
\begin{equation}\displaystyle
\langle \frac{\partial}{\partial \overline{z}}f,g\rangle_{\mathcal{SF}(\mathbb{C})}=\langle f,M_{\overline{z}}g\rangle_{\mathcal{SF}(\mathbb{C})}.
\end{equation}
\end{theorem}
\begin{proof}
Let $\displaystyle f=\sum_{k=0}^{\infty}\overline{z}^kf_k$ and $\displaystyle g=\sum_{k=0}^{\infty}\overline{z}^kg_k$ in $\mathcal{SF}(\mathbb{C})$ that belongs to the domains of the creation and annihilation operators. Firstly, we note that we have $$\displaystyle M_z(g)=\sum_{k=0}^{\infty}\overline{z}^kM_z(g_k).$$ Then, it follows that
\[
    \begin{split}
   \displaystyle  \langle \frac{\partial}{\partial z}f,g\rangle_{\mathcal{SF}(\mathbb{C})}
      &= \sum_{k=0}^{\infty}k!\langle \frac{\partial}{\partial z} f_k,g_k \rangle_{\mathcal{F}(\mathbb{C})}\\
        &= \sum_{k=0}^{\infty}k!\langle f_k,(\frac{\partial}{\partial z})^*g_k \rangle_{\mathcal{F}(\mathbb{C})}
        \\
        &= \sum_{k=0}^{\infty}k!\langle f_k,M_z(g_k) \rangle_{\mathcal{F}(\mathbb{C})}
         \\
        &=\langle f,M_{z}(g)\rangle_{\mathcal{SF}(\mathbb{C})}.
          \end{split}
   \]
   Moreover, since $$\displaystyle \frac{\partial}{\partial\overline{z}}(f)(z)=\sum_{h=0}^{\infty}(h+1)\overline{z}^hf_{h+1},$$
   and
   $$\displaystyle M_{\overline{z}}(g)(z)=\sum_{h=1}^{\infty}\overline{z}^hg_{h-1},$$
   it follows that
   \[
    \begin{split}
   \displaystyle  \langle \frac{\partial}{\partial\overline{z}}f,g\rangle_{\mathcal{SF}(\mathbb{C})}    &= \sum_{k=0}^{\infty}k!\langle (k+1)f_{k+1},g_k \rangle_{\mathcal{F}(\mathbb{C})}
        \\
        &= \sum_{k=0}^{\infty}(k+1)!\langle f_{k+1},g_k \rangle_{\mathcal{F}(\mathbb{C})}
         \\
           &= \sum_{k=1}^{\infty}k!\langle f_{k},g_{k-1} \rangle_{\mathcal{F}(\mathbb{C})}
         \\
        &=\langle f,M_{\overline{z}}(g)\rangle_{\mathcal{SF}(\mathbb{C})}.
          \end{split}
   \]
\end{proof}
\begin{remark}
We shall see later that that we have  $$\mathcal{F}(\mathbb{C})\subset \mathcal{SF}(\mathbb{C})\subset L^2(\mathbb{C},d\mu_\beta), \quad \beta>2.$$
The previous inclusions are strict. It is well known that the classical Fock space $\mathcal{F}(\mathbb{C})$ is the only space of entire functions on which the creation and annihilation operators are adjoints of each others and satisfy the classical commutation rules. Of course, this is not true anymore on $L^2(\mathbb{C},d\mu)$, see Proposition 7.2 in \cite{Shige}. However, the previous theorem shows that the result still holds in the subspace $\mathcal{SF}(\mathbb{C})$ of the space of polyanalytic functions of infinite order.
\end{remark}

\begin{definition}\label{RinfinityL}
Let $f,g\in\mathcal{SF}(\mathbb{C})$ and let $f(z)=\displaystyle \sum_{n=0}^{\infty}\bar{z}^nf_n(z)$ and $g(z)=\displaystyle \sum_{m=0}^{\infty}z^mg_m(\bar{z})$. Then, we define two backward shift operators $R_\infty$ and $L_\infty$ with respect to the variables $z$ and $\bar{z}$ respectively given by
\begin{equation}
R_{\infty}(f)(z,\bar{z})=\displaystyle \sum_{n=0}^{\infty}\bar{z}^nR_0(f_n)(z)=\frac{1}{z}\left(f(z,\bar{z})-\sum_{n=0}^{\infty}\bar{z}^{n}f_n(0)\right), \quad z\in\mathbb{C}
\end{equation}
and
\begin{equation}
L_{\infty}(g)(z,\bar{z})=\displaystyle \frac{1}{\bar{z}}\left(g(z,\bar{z})-\sum_{m=0}^{\infty}z^{m}g_m(0)\right), \quad z\in\mathbb{C}.
\end{equation}
\end{definition}
It turns out that both the backward shift operators $R_\infty$ and $L_\infty$ define two contractions on the polyanalytic Fock space of infinite order $\mathcal{SF}(\mathbb{C})$. Indeed, the following result holds:
\begin{proposition}
For any $f\in\mathcal{SF}(\mathbb{C}),$ we have
\begin{equation}
||R_{\infty}(f)||_{\mathcal{SF}(\mathbb{C})}^2\leq ||f||_{\mathcal{SF}(\mathbb{C})}^2-\sum_{n=0}^{\infty}n!|f_n(0)|^2 \leq||f||_{\mathcal{SF}(\mathbb{C})}^2.
\end{equation}
and \begin{equation}
||L_{\infty}(g)||_{\mathcal{SF}(\mathbb{C})}^2\leq ||g||_{\mathcal{SF}(\mathbb{C})}^2-\sum_{n=0}^{\infty}n!|g_n(0)|^2 \leq||g||_{\mathcal{SF}(\mathbb{C})}^2.
\end{equation}
\end{proposition}
\begin{proof}
We will prove the result for $R_{\infty}$. Indeed, if we consider
$$f(z)=\displaystyle \sum_{n=0}^{\infty}\bar{z}^n f_n(z), \quad f_n\in\mathcal{F}(\mathbb{C}),$$
we have $$R_{\infty}(f)(z) \displaystyle =\sum_{n=0}^{\infty}\bar{z}^nR_0(f_{n})(z).$$
Thus, using the fact $R_0$ is a contraction on the Fock space (see \cite{ACK}) we deduce
  \[
    \begin{split}
   \displaystyle ||R_{\infty}(f)||_{\mathcal{SF}(\mathbb{C})}^2    &= \sum_{n=0}^{\infty}n!||R_0(f_{n})||_{\mathcal{F}(\mathbb{C})}^2
        \\
        &\leq  \sum_{n=0}^{\infty}n!(||f_{n}||_{\mathcal{F}(\mathbb{C})}^2-|f_n(0)|^2)
         \\
           &=||f||_{\mathcal{SF}(\mathbb{C})}^2-\sum_{n=0}^{\infty}n!|f_n(0)|^2\leq ||f||_{\mathcal{SF}(\mathbb{C})}^2
         \\
        &
          \end{split}
   \]
   and this ends the proof.  We note that the argument follows in a similar way for $L_{\infty}$ using the fact that we can write $f$ also in the form
$$f(z)=\displaystyle\sum_{m=0}^{\infty}z^mf_m(\overline{z}).$$
\end{proof}

Now, we consider other two operators:

\begin{definition}
Let $f,g\in\mathcal{SF}(\mathbb{C})$ and let $f(z)=\displaystyle \sum_{n=0}^{\infty}\bar{z}^nf_n(z)$ and $g(z)=\displaystyle \sum_{m=0}^\infty z^mg_m(\bar{z})$. Then, we define two operators with respect to the variables $z$ and $\bar{z}$ respectively which are given by
\begin{equation}
I_{\infty}(f)(z,\bar{z})=\displaystyle \sum_{n=0}^{\infty}\bar{z}^nI(f_n)(z), \quad z\in\mathbb{C}
\end{equation}
and
\begin{equation}
J_{\infty}(g)(z,\bar{z})=\displaystyle \sum_{m=0}^{\infty}z^mJ(g_m)(\bar{z}), \quad z\in\mathbb{C}.
\end{equation}
\end{definition}
We point out that the operator $I$ is the integration operator considered in \cite{ACK}, while $J$ is the integration with respect to the conjugate variable.

\begin{remark}
It holds that
\begin{equation}
I_{\infty}(\bar{z}^nz^m)(z)=\displaystyle \frac{\bar{z}^{n}z^{m+1}}{m+1} , \quad z\in\mathbb{C}
\end{equation}
and
\begin{equation}
J_{\infty}(\bar{z}^nz^m)(z)= \displaystyle \frac{\bar{z}^{n+1}z^m}{n+1} , \quad z\in\mathbb{C}.
\end{equation}
\end{remark}
As a consequence, we have the following result:
\begin{theorem}
The adjoints of $R_\infty$ and $L_\infty$ satisfy
\begin{equation}
R_{\infty}^{*}=I_\infty
\end{equation}
and
\begin{equation}
L_{\infty}^{*}=J_\infty.
\end{equation}
\end{theorem}
\begin{proof}
Let $f,g\in\mathcal{SF}(\mathbb{C})$; we will prove that

\begin{equation}
\langle I_{\infty}(f),g \rangle_{\mathcal{SF}(\mathbb{C})}=\langle f, R_{\infty}(g) \rangle_{\mathcal{SF}(\mathbb{C})}.
\end{equation}
Indeed, we write $f(z)=\displaystyle \sum_{n=0}^{\infty}\overline{z}^n f_n(z)$ and $g(w)=\displaystyle \sum_{n=0}^{\infty}\overline{z}^n g_n(z)$ with $f_n,g_n\in\mathcal{F}(\mathbb{C})$ for any $n\geq 0$. We have

$$I_{\infty}(f)(z)=\displaystyle \sum_{n=0}^{\infty} \bar{z}^nI(f_n)(z) $$
and
$$R_{\infty}(g)(z)=\sum_{m=0}^\infty \bar{z}^m R_0(g_m)(z).$$

Then, we use the scalar product on $\mathcal{SF}(\mathbb{C})$ and apply the result on the classical backward shift operator, see \cite{ACK}.  Therefore, it follows that
   \[
    \begin{split}
   \displaystyle  \langle I_\infty(f),g\rangle_{\mathcal{SF}(\mathbb{C})}    &= \sum_{k=0}^{\infty}k!\langle I(f_{k}),g_k \rangle_{\mathcal{F}(\mathbb{C})}
        \\
        &= \sum_{k=0}^{\infty}k!\langle f_{k},R_0(g_k) \rangle_{\mathcal{F}(\mathbb{C})}
         \\
           &= \langle f,R_\infty(g) \rangle_{\mathcal{SF}(\mathbb{C})}.
         \\
        &
          \end{split}
   \]
   The second part of the statement can be proved in a similar way.
\end{proof}

\section{A kernel function relating polyanalytic Fock spaces of finite and infinite order}
\setcounter{equation}{0}
In this section we study how the polyanalytic Fock spaces of finite and infinite order are related between them.
We denote by $\mathcal{F}_n(\mathbb{C})$ the classical polyanalytic Fock space whose kernel is given by the formula \eqref{Kn}.
The relation between the kernels $K$ and $(K_n)_{n\geq 1}$ is described in the next result.
\begin{proposition}[kernel formula]
For any $z,w\in\mathbb{C}$ we set $$G(z,w)=e^{z\overline{w}-(|z|^2+|w|^2)}.$$
Then, it holds that
\begin{equation}
\displaystyle  \sum_{n=1}^{\infty}  \frac{K_n(z,w)}{2^{n+1}}=G(z,w)K(z,w), \quad \text{ for any } z,w\in\mathbb{C}.
\end{equation}

\end{proposition}
\begin{proof}
We note that the polyanalytic Fock kernels given by \eqref{Kn} can be written in terms of the generalized Laguerre polynomials as follows
\begin{equation}\label{LKn}
K_n(z,w)=e^{z\overline{w}}L_{n-1}^{1}(|z-w|^2), \quad \text{ for any } z,w\in\mathbb{C}.
\end{equation}
Then, taking the series \eqref{LKn} we obtain

\begin{equation}\label{sumKn}
\displaystyle \sum_{n=1}^{\infty}\frac{K_n(z,w)}{2^{n-1}}=e^{z\overline{w}}\sum_{n=1}^{\infty}\frac{L_{n-1}^{1}(|z-w|^2)}{2^{n-1}}=e^{z\overline{w}}\sum_{n=0}^{\infty}\frac{L_{n}^{1}(|z-w|^2)}{2^{n}}.
\end{equation}
Morever, we note that for any $a,\alpha>0$ we have the following expansion, see \cite[Example 2, pp 89]{L}

\begin{equation}\label{ax}
e^{-ax}=(a+1)^{-(\alpha+1)}\sum_{n=0}^{\infty}\left(\frac{a}{a+1}\right)^{n}L_n^\alpha(x),\quad x\geq 0.
\end{equation}
In particular, inserting  $\alpha=a=1$ and $x=|z-w|^2$ in \eqref{ax} we obtain
$$e^{-|z-w|^2}=\frac{1}{2^2}\sum_{n=0}^{\infty}\frac{L^1_n(|z-w|^2)}{2^n}, \quad z,w\in\mathbb{C}.$$
Hence, with some computations involving \eqref{sumKn} we obtain that

$$\displaystyle  \sum_{n=1}^{\infty}  \frac{K_n(z,w)}{2^{n+1}}=G(z,w)K(z,w),$$
 \text{ for any } $z,w\in\mathbb{C}$ where $G(z,w)=e^{z\overline{w}-(|z|^2+|w|^2)}.$
\end{proof}

\begin{remark}
We observe that the classical creation and annihilation operators are adjoint of each others on the polyanalytic Fock space of infinite order $\mathcal{SF}(\mathbb{C})$, see Theorem \ref{Adj}.
\end{remark}
\section{A Segal-Bargmann type transform and related operators}
\setcounter{equation}{0}
In this section, we deal with a Segal-Bargmann type transform related to the polyanalytic Fock spaces of infinite order. We discuss also some related operators.

Let ($\psi_n(x))_{n\geq 0}$ denote the normalized Hermite functions and consider the Segal-Bargmann kernel $A(z,x)$ which is given by \begin{equation}
A(z,x):=\displaystyle \sum_{n=0}^{\infty}\frac{z^n}{\sqrt{n!}}\psi_n(x)=e^{-\frac{1}{2}(z^2+x^2)+\sqrt{2}zx}, \quad \text{ for any } (z,x)\in\mathbb{C}\times \mathbb{R}.
\end{equation}
For any $z\in\mathbb{C}$  fixed we use also the notation $A_z(x)=A(z,x)$ for all $x\in\mathbb{R}.$ The kernel \eqref{newFock-kernel} can be factorized as follows:
\begin{theorem} \label{kerFac}
For any $(z,w)\in\mathbb{C}^2$, we have

\begin{equation}
K(z,w)= \langle A_z\otimes A_{\bar{z}},A_w \otimes A_{\bar{w}}\rangle_{L^2(\mathbb{R}^2)}.
\end{equation}
\end{theorem}
\begin{proof}
The proof is based on computations using Fubini's theorem combined with the following well-known fact  $$\langle A_z, A_w\rangle_{L^2(\mathbb{R})}=e^{z\bar{w}}.$$
Indeed, for $z,w\in\mathbb{C}$ we have the explicit computations

 \[
    \begin{split}
   \displaystyle   \langle A_z\otimes A_{\bar{z}},A_w \otimes A_{\bar{w}}\rangle_{L^2(\mathbb{R}^2)} &= \int_{\mathbb{R}^2} (A_z\otimes A_{\bar{z}})(x,y)\overline{(A_w\otimes A_{\bar{w}})(x,y)} dxdy\\
      &=\int_{\mathbb{R}^2}A_z(x)A_{\overline{z}}(y) \overline{A_w(x)}\textbf{  }\overline{A_{\overline{w}}(y)}dxdy\\
        &= \left(\int_{\mathbb{R}}A_z(x) \overline{A_w(x)}dx \right) \left(\int_{\mathbb{R}}\overline{A_{z}(y)} A_w(y) dy\right)
        \\
        &= \langle A_z, A_w\rangle_{L^2(\mathbb{R})} \langle A_w, A_z\rangle_{L^2(\mathbb{R})}
        \\
        &= e^{z\bar{w}} e^{w\bar{z}}
         \\
        &=K(z,w).
          \end{split}
   \]

\end{proof}

\begin{definition}\label{Ttran}
For a given $\varphi\in L^2(\mathbb{R}^2)$, we define the so-called first Segal-Bargmann type transform by
\begin{equation} \displaystyle
T(\varphi)(z,\bar{z})=\langle \varphi, \overline{A_z\otimes A_{\bar{z}}}\rangle_{L^2(\mathbb{R}^2)}= \int_{\mathbb{R}^2} A_z(x)A_{\bar{z}}(y)\varphi(x,y)dxdy.
\end{equation}
\end{definition}
Then, as a consequence, we can write the kernel function $K(z,w)$
as a function in the range of the transform $T$ thanks to the following

\begin{proposition}
For a fixed $w\in\mathbb{C}$, we set $\varphi_w(t_1,t_2)=(A_{\overline{w}}\otimes A_w)(t_1,t_2)$, with $t_1,t_2\in\mathbb{R}$. Then

$$K(z,w)=T(\varphi_w)(z,\overline{z}), \quad \text{ for any } z\in\mathbb{C}.$$

\end{proposition}
\begin{proof}
This result can be obtained as a direct application of Theorem \ref{kerFac} and Definition \ref{Ttran} taking into account that
$$\overline{A_{\overline{w}}\otimes A_w}=A_w\otimes A_{\overline{w}}, \quad w\in\mathbb{C}.$$

 Indeed, for any $z,w\in\mathbb{C}$ we have
 \[
    \begin{split}
   \displaystyle  T(\varphi_w)(z,\overline{z})&= \langle \varphi_w, \overline{A_z\otimes A_{\bar{z}}}\rangle_{L^2(\mathbb{R}^2)}\\
      &= \langle A_{\overline{w}}\otimes A_{w}, \overline{A_z\otimes A_{\bar{z}}}\rangle_{L^2(\mathbb{R}^2)}\\
        &= \int_{\mathbb{R}^2}(A_{\overline{w}}\otimes A_{w})(t_1,t_2)(A_z\otimes A_{\bar{z}})(t_1,t_2)dt_1dt_2
        \\
        &= \langle A_z\otimes A_{\bar{z}}, \overline{A_{\overline{w}} \otimes A_{w}}\rangle_{L^2(\mathbb{R}^2)}
        \\
        &=  \langle A_z\otimes A_{\bar{z}},A_w \otimes A_{\bar{w}}\rangle_{L^2(\mathbb{R}^2)}
         \\
        &=K(z,w).
          \end{split}
   \]

\end{proof}

\begin{remark}
As a consequence of the previous result, we observe that for any fixed $w\in\mathbb{C}$,  we have

$$T(\varphi_w)(z)=K(z,w)=K_w(z), \quad z\in\mathbb{C}.$$

Then, $T^*=T^{-1}$ since $T$ is a unitary operator; moreover, for any $w\in\mathbb{C}$ 

$$T^{-1}(K_w)(t_1,t_2)=T^{*}(K_w)(t_1,t_2)=(A_{\overline{w}}\otimes A_w)(t_1,t_2), \quad  \text{ for all } (t_1,t_2)\in\mathbb{R}^2.$$
\end{remark}

As a first example, we consider the family of functions given by $$\psi_{m,n}(x,y):=(\psi_m\otimes \psi_n)(x,y)=\psi_m(x)\psi_n(y), \text{ for any } m,n\geq 0.$$ We have
\begin{proposition}\label{Action1}
For every $z\in\mathbb{C}$  we have
\begin{equation}
T(\psi_{m,n})(z,\bar{z})=\frac{z^m\overline{z^n}}{\sqrt{m!n!}}=\phi_{m,n}(z,\overline{z}), \quad m,n=0,1,...
\end{equation}
and
\begin{equation}
\Delta_z \phi_{p,q}(z,\overline{z})=4\sqrt{pq}\phi_{p-1,q-1}(z,\overline{z}), \quad p,q=1,2,....
\end{equation}
\end{proposition}
\begin{proof}
Recalling that $A_{\overline{z}}(y)=\overline{A_z(y)}$ for any $y\in\mathbb{R}$, using the Fubini's theorem we have:   \[
    \begin{split}
   \displaystyle T(\psi_{(m,n)})(z,\bar{z}) &= \int_{\mathbb{R}^2} (A_z\otimes A_{\bar{z}})(x,y)\psi_{m,n}(x,y)dxdy\\
      &=\int_{\mathbb{R}^2}A_z(x)A_{\overline{z}}(y) \psi_{m}(x)\psi_n(y) dxdy\\
        &= \left(\int_{\mathbb{R}}A_z(x) \psi_m(x) dx \right) \left(\int_{\mathbb{R}}\overline{A_{z}(y)}\psi_n(y) dy\right)
        \\
        &= B(\psi_m)(z)\overline{B(\psi_n)(z)}
        \\
        &= \frac{z^m\overline{z}^n}{\sqrt{m!n!}}
         \\
        &=\phi_{m,n}(z,\overline{z}).
          \end{split}
   \]

   Now, using the fact that $\displaystyle \Delta_z=4\frac{\partial^2}{\partial z  \partial \overline{z}}$ we get

   $$\displaystyle \Delta_z(\phi_{p,q})(z,\overline{z})=4\frac{pq}{\sqrt{p!q!}}z^{p-1}\overline{z}^{q-1}=4\sqrt{pq}\phi_{p-1,q-1}(z,\overline{z}).$$

\end{proof}
\begin{example}
Let $$f(z,\overline{z})=T(\psi_{n,m})(z,\overline{z}), \text{ for any } z\in\mathbb{C}.$$
 Then $f\in\mathcal{SF}(\mathbb{C})$, moreover we have $$||f||_{\mathcal{SF}(\mathbb{C})}=||T(\psi_{n,m})||=1=||\psi_{n,m}||_{L^2(\mathbb{R}^2)}.$$
 \end{example}
\begin{theorem}
The first Segal-Bargmann type transform $T$ defines an isometric isomorphism from $L^2(\mathbb{R}^2)$ onto $\mathcal{SF}(\mathbb{C})$.
\end{theorem}
\begin{proof}
The normalized Hermite functions $(\psi_{m,n})_{m,n\geq 0}$ form an orthonormal basis of $L^2(\mathbb{R}^2)$, thus for any $\varphi\in L^2(\mathbb{R}^2,\mathbb{C}),$  there exist unique coefficients $(\beta_{m,n})_{m,n\geq 0}$ in $\mathbb{C}$ such that
$$\varphi(x,y)=\displaystyle \sum_{m,n=0}^{\infty} \psi_{m,n}(x,y)\beta_{m,n}, \quad \text{ and }\quad ||\varphi||_{L^2(\mathbb{R}^2)}^{2}=\sum_{m,n=0}^{\infty}|\beta_{m,n}|^2<\infty. $$
Therefore, inserting $\varphi$ in the definition of the transform $T$  and using some standard arguments we have

 \[
    \begin{split}
   \displaystyle T(\varphi)(z,\bar{z}) &= \int_{\mathbb{R}^2} A_z(x) A_{\bar{z}}(y)\left(\sum_{m,n=0}^{\infty}\psi_{m,n}(x,y)\beta_{m,n}\right)dxdy\\
      &=\sum_{m,n=0}^{\infty}\left(\int_{\mathbb{R}^2}A_z(x)A_{\overline{z}}(y) \psi_{m,m}(x,y) dxdy\right) \beta_{m,n}\\
        &= \sum_{m,n=0}^{\infty}T(\psi_{m,n})(z,\overline{z})\beta_{m,n}.
            \end{split}
   \]

Applying the first part of Proposition \ref{Action1} we obtain
 \[
    \begin{split}
   \displaystyle T(\varphi)(z,\bar{z}) &= \sum_{m,n=0}^{\infty}\phi_{m,n}(z,\overline{z})\beta_{m,n}
        \\
        &= \sum_{m,n=0}^{\infty} \frac{z^m\overline{z}^n}{\sqrt{m!n!}}\beta_{m,n}
         \\
        &= \sum_{n=0}^{\infty}\overline{z}^n\left(\sum_{m=0}^{\infty}\frac{z^m}{\sqrt{m!n!}}\beta_{m,n}\right).
          \end{split}
   \]
Then, setting $\displaystyle f_n(z)=\frac{1}{\sqrt{n!}}\sum_{m=0}^{\infty}\frac{z^m}{\sqrt{m!}}\beta_{m,n}$ for any $n\geq 0 $ we have

$$T(\varphi)(z,\overline{z})=\displaystyle \sum_{n=0}^{\infty}\overline{z}^nf_n(z), \quad z\in\mathbb{C}.$$
We observe that $f_n$ are entire functions. Moreover, it is immediate to see that
$$\displaystyle ||f_n||^{2}_{\mathcal{F}(\mathbb{C})}=\frac{1}{n!}\sum_{m=0}^{\infty}|\beta_{m,n}|^2<\infty.$$
Hence, we deduce
  \[
    \begin{split}
   \displaystyle ||T(\varphi)||_{\mathcal{SF}(\mathbb{C})}^{2} &= \sum_{n=0}^{\infty} n!||f_n||^{2}_{\mathcal{F}(\mathbb{C})}\\
      &=\sum_{n=0}^{\infty}\sum_{m=0}^{\infty}|\beta_{m,n}|^2\\
        &=||\varphi||_{L^2(\mathbb{R})}^{2}.
        \\
        &
          \end{split}
   \]
   This proves that the transform $T$ defines an isometric operator from $L^2(\mathbb{R})$ into the polyanalytic Fock space of infinite order $\mathcal{SF}(\mathbb{C}).$ On the other hand, we note that using Proposition \ref{Action1} we have

  $$T(\psi_{m,n})=\frac{z^m\overline{z}^n}{\sqrt{m!n!}}, \quad m,n\geq 0.$$
  This allows to justify that if $f\in\mathcal{SF}(\mathbb{C})$ there exists $g\in L^2(\mathbb{R})$ such that $f=T(g)$. In particular, this shows that $T$ is also surjective which ends the proof.
\end{proof}
\begin{remark}
We know that $T$ defines an isometric and surjective transform from $L^2(\mathbb{R}^2)$ onto $\mathcal{SF}(\mathbb{C})$. This shows that $T$ is invertible and the inverse of $T$ exists. However, we do not know how to explicitly calculate this inverse because of the lack of a geometric description of the space  $\mathcal{SF}(\mathbb{C})$.
\end{remark}
We recall the position operators given by $$X\varphi(x,y)=x\varphi(x,y)$$  and $$Y\varphi(x,y)=y\varphi(x,y).$$

We denote by $\mathcal{D}(X)=\left\lbrace{\varphi\in L^2(\mathbb{R}), X(\varphi)\in L^2(\mathbb{R})}\right\rbrace$ the domain of the position operator $X$.
\begin{proposition}\label{PositionX}
 The following relations hold:
$$\displaystyle T^{-1}\left(\frac{\partial}{\partial z}+M_z\right)T=\sqrt{2}X, \quad \text{ on  } \mathcal{D}(X)$$

and $$\displaystyle T^{-1}\left(\frac{\partial }{\partial \overline{z}}+M_{\overline{z}}\right)T=\sqrt{2}Y, \quad \text{ on } \mathcal{D}(Y).$$
\end{proposition}
\begin{proof}
We will make the calculations for the operator $X$ and for $Y$ it can be done in a similar way. As in the classical case, for $\varphi\in \mathcal{D}(X)$ we have
$$\displaystyle \frac{\partial}{\partial z}[T(\varphi)](z)=\displaystyle \int_{\mathbb{R}^2}\frac{\partial}{\partial z}(A_z(x))A_{\overline{z}}(y)\varphi(x,y)dxdy$$.

It is easy to check that

$$\displaystyle \frac{\partial}{\partial z} (A(z,x))=(-z+\sqrt{2}x)A(z,x).$$
Thus, we obtain

 \[
    \begin{split}
   \displaystyle  \frac{\partial}{\partial z} T(\varphi)(z)  &= \int_{\mathbb{R}^2}(-z+\sqrt{2}x)A_z(x)A_{\overline{z}}(y)\varphi(x,y)dxdy
        \\
        &= -z \int_{\mathbb{R}^2}A_z(x)A_{\overline{z}}(y)\varphi(x,y)dxdy+ \sqrt{2}\int_{\mathbb{R}^2}A_z(x)A_{\overline{z}}(y)x\varphi(x,y)dxdy
         \\
           &=-M_zT(\varphi)(z)+\sqrt{2}T(X \varphi)(z).
         \\
        &
          \end{split}
   \]

   Therefore, it follows that for any $\varphi \in \mathcal{D}(X)$ we have
   $$\displaystyle \frac{\partial}{\partial z} T(\varphi)+M_zT(\varphi)=\sqrt{2}TX(\varphi).$$
   Hence, we obtain
   $$\left(\frac{\partial}{\partial z}+M_z \right)T=\sqrt{2}TX.$$
   Finally, this leads to
   $$\displaystyle T^{-1}\left(\frac{\partial}{\partial z}+M_z\right)T=\sqrt{2}X, \quad \text{ on  } \mathcal{D}(X).$$

\end{proof}
\begin{proposition}
We have
 $$\frac{1}{2}T^{-1}\left(\frac{1}{4}\Delta_z+z\frac{\partial}{\partial \overline{z}}+\overline{z}\frac{\partial}{\partial z}+M_{|z|^2}\right)T=YX, \text{ on } \mathcal{D}(YX).$$

\end{proposition}
\begin{proof}
We observe that Proposition \ref{PositionX} yields
$$\displaystyle\frac{1}{2}T^{-1}\left(\frac{\partial}{\partial \overline{z}}+M_{\overline{z}}\right)\left(\frac{\partial}{\partial z}+M_z\right)T=YX.$$

Then, the result holds since we have $$\left(\frac{\partial}{\partial\overline{z}}+M_{\overline{z}}\right)\left(\frac{\partial}{\partial z}+M_z\right)=\frac{1}{4}\Delta_z+z\frac{\partial}{\partial \overline{z}}+\overline{z}\frac{\partial}{\partial z}+M_{|z|^2}.$$

\end{proof}
Let $\mathcal{D}(M_z)=\left\lbrace{f\in \mathcal{SF}(\mathbb{C}), M_z(f)\in \mathcal{SF}(\mathbb{C})}\right\rbrace$ be the domain of the creation operator $M_z$.
\begin{proposition} \label{CreationM}
The following relations hold
 $$T\left(X-\frac{\partial}{\partial x}\right)T^{-1}=\sqrt{2}M_z, \quad \text{ on }\mathcal{D}(X)\cap \mathcal{D}(\frac{\partial}{\partial x})$$
 and
 $$T\left(Y-\frac{\partial}{\partial y}\right)T^{-1}=\sqrt{2}M_{\overline{z}}, \quad \text{ on }\mathcal{D}(Y)\cap \mathcal{D}(\frac{\partial}{\partial y}).$$
\end{proposition}
\begin{proof}
We will prove the first statement of this result; the second one follows with similar arguments. Indeed, let $\displaystyle \varphi\in \mathcal{D}(X)\cap \mathcal{D}(\frac{\partial}{\partial x})$ we have
$$T(X-\frac{\partial}{\partial x})\varphi (z)=T(X\varphi)(z)-T(\frac{\partial}{\partial x}\varphi)(z).$$

By definition we have
$$\displaystyle T(X\varphi)(z)=\displaystyle \int_{\mathbb{R}^2} A_z(x)A_{\bar{z}}(y)x\varphi(x,y)dxdy,$$

so that
$$\displaystyle  T(\frac{\partial}{\partial x}\varphi)(z)= \int_{\mathbb{R}^2} A_z(x)A_{\bar{z}}(y)\frac{\partial}{\partial x}\varphi(x,y)dxdy=- \int_{\mathbb{R}^2} \frac{\partial}{\partial x}(A_z(x))A_{\bar{z}}(y)\varphi(x,y)dxdy,$$

since $$\displaystyle \frac{\partial}{\partial x}(A_z(x))=(-x+\sqrt{2}z)A_z(x),$$

we are lead to
$$T\left(X-\frac{\partial}{\partial x}\right)\varphi(z)=\sqrt{2}zT(\varphi)(z).$$
 Hence, we obtain

 $$T\left(X-\frac{\partial}{\partial x}\right)T^{-1}=\sqrt{2}M_z, \quad \text{ on }\mathcal{D}(X)\cap \mathcal{D}(\frac{\partial}{\partial x}).$$
\end{proof}
As a consequence of the previous result we can prove the following
\begin{corollary}
We have
 $$T\left(\frac{\partial^2}{\partial x\partial y}+XY-X\frac{\partial}{\partial y}-Y\frac{\partial}{\partial x}\right)T^{-1}=2M_{|z|^2}.$$
\end{corollary}
\begin{proof}
Indeed, we just need to apply Proposition \ref{CreationM} combined with the relation $$\left(X-\frac{\partial}{\partial x}\right)\left(Y-\frac{\partial}{\partial y}\right)=\frac{\partial^2}{\partial x \partial y}+XY-X\frac{\partial}{\partial y}-Y\frac{\partial}{\partial x}.$$
\end{proof}
\section{A Berezin transform and related operators}
\setcounter{equation}{0}
We now use the kernel function \eqref{newFock-kernel} to study a Berezin integral transform and develop further results on it.
\begin{definition}
Let $f:\mathbb{C}\longrightarrow \mathbb{C}$ and let $d\mu(w)=\frac{1}{\pi}e^{-|w|^2}dA(w)$ be the Gaussian measure. Then, we consider the following integral transform
\begin{equation}\label{Noperator}
\displaystyle \mathcal{B}(f)(z)=\int_{\mathbb{C}} e^{-|z|^2}K(z,w)f(w)d\mu(w), \quad z\in\mathbb{C},
\end{equation}
when the integral exists.
\end{definition}
\begin{remark}
We observe that $$e^{-|z|^2}K(z,w)e^{-|w|^2}=e^{-|z-w|^2},\quad z,w\in \mathbb{C}.$$
Thus, it turns out that the integral transform given by \eqref{Noperator} coincides with the so-called Berezin transform considered in \cite[p. 101]{Zhu}. However, since $\mathcal{B}$ can be expressed in terms of the kernel function $K(z,w)$ we can use various properties of this kernel in order to develop further results. A similar transform in the case of two complex variables was considered in \cite{BG2019}. It is important to note that the Berezin transform was first introduced by Berezin in \cite{Berezin} as a general concept of quantization.
\end{remark}
We start first by observing that the Berezin transform $\mathcal{B}$ fixes  all the complex monomials $z^n$, $n\in\mathbb N$ which form an orthogonal basis of the classical Fock $\mathcal{F}(\mathbb{C})$:
\begin{proposition}\label{Bzn}
For $n=0,1,...$ it holds that
\begin{equation}
\mathcal{B}(z^n)=z^n,\quad \forall z\in\mathbb{C}.
\end{equation}

\end{proposition}

\begin{proof}
We set $f_n(z)=z^n$ with $n=0,1,...$.
Then, for every $\alpha\in\mathbb{C}$ we can write
\begin{equation}
\displaystyle \sum_{n=0}^{\infty}\alpha^n \frac{\mathcal{N}(f_n)(z)}{n!}=\frac{1}{\pi}e^{-|z|^2}\int_{\mathbb{C}}K(z,w)e^{\alpha w}e^{-|w|^2}dA(w).
\end{equation}
Let us set $w=t_1+it_2$ and let us replace $K(z,w)$ by its expression to obtain
$$\displaystyle \sum_{n=0}^{\infty}\alpha^n \frac{\mathcal{N}(f_n)(z)}{n!}=\frac{1}{\pi}e^{-|z|^2}\int_{\mathbb{R}^2}e^{z(t_1-it_2)+(t_1+it_2)\overline{z}+\alpha(t_1+it_2)}e^{-(t_1^2+t_2^2)}dt_1dt_2.$$
Therefore, thanks to Fubini's theorem we have
$$\displaystyle \sum_{n=0}^{\infty} \alpha^n \frac{\mathcal{B}(f_n)(z)}{n!}=\frac{1}{\pi}e^{-|z|^2}\left(\int_{\mathbb{R}}e^{-t_1^2+t_1(z+\overline{z}+\alpha)}dt_1\right)\left(\int_{\mathbb{R}}e^{-t_2^2+t_2i(\overline{z}-z+\alpha)}dt_2\right).$$
Now, we recall the classical Gaussian integral

\begin{equation}
\displaystyle \int_{\mathbb{R}}e^{-at^2+bt}dt=\sqrt{\frac{\pi}{a}}e^{\frac{b^2}{4a}},\qquad a>0, b\in\mathbb{C}.
\end{equation}
Thus, setting $b_1(z,\overline{z})=z+\overline{z}+\alpha$ and $b_2(z,\overline{z})=i(\overline{z}-z+\alpha)$ we obtain
$$\displaystyle \sum_{n=0}^{\infty} \alpha^n\frac{\mathcal{B}(f_n)(z)}{n!}=\frac{1}{\pi}e^{-|z|^2}\pi e^{\frac{b_1^2(z,\overline{z})}{4}}e^{\frac{b_2^2(z,\overline{z})}{4}}.$$

Therefore, we have $$\displaystyle \sum_{n=0}^{\infty} \alpha^n\frac{\mathcal{B}(f_n)(z)}{n!}=e^{-|z|^2}e^{\frac{b^2_1(z,\overline{z})+b^2_2(z,\overline{z})}{4}}$$

We have $b_1^2(z,\overline{z})=z^2+\overline{z}^2+2|z|^2+\alpha^2+2\alpha z+2\alpha \overline{z}$ and $b_2^2(z,\overline{z})=-(\overline{z}^2+z^2-2|z|^2+\alpha^2+2\alpha \overline{z}-2\alpha z)$, which leads to

$$b^2_1(z,\overline{z})+b^2_2(z,\overline{z})=4(|z|^2+\alpha z).$$
Hence
$$\displaystyle \sum_{n=0}^{\infty} \alpha^n\frac{\mathcal{B}(f_n)(z)}{n!}=e^{\alpha z}=\sum_{n=0}^{\infty}\alpha^n \frac{z^n}{n!}, \quad \forall \alpha \in\mathbb{C}.$$

Finally, we identify the coefficients with respect to the variable $\alpha$ and get
$$\mathcal{B}(f_n)(z)=z^n, \quad  n=0,1,...$$
\end{proof}
\begin{remark}
It is impotant to note that in \cite{Zhu} it was proved that a function $f\in L^\infty(\mathbb{C})$ is a fixed point for the Berezin transform $\mathcal{B}$ if and only if $f$ is constant. Then, Proposition \ref{Bzn} shows that monomials are fixed by the Berezin transform. In particular, this suggests to consider a more general fixed point problem: to find when $\mathcal{B}(f)=f$ for $f$ in some suitable function space.
\end{remark}
\begin{remark}
We observe that
\begin{equation}
\displaystyle \mathcal{B}(f)(z)=\frac{1}{\pi}\int_{\mathbb{C}} e^{\overline{w}z+w\overline{z}-|z|^2}f(w)e^{-|w|^2}dA(w).
\end{equation}
If $z=0$, we have
$$\displaystyle \mathcal{B}(f)(0)=\frac{1}{\pi}\int_{\mathbb{C}} f(w)e^{-|w|^2}dA(w).$$
In particular, using the Cauchy-Schwarz inequality we obtain
$$|\mathcal{B}(f)(0)|\leq ||f||_{L^2(\mathbb{C},\mu)},$$
and by taking the function $f=1$, we get
$$\mathcal{B}(1)(z)=\displaystyle e^{-|z|^2} \frac{1}{\pi}\int_{\mathbb{C}} K(z,w)e^{-|w|^2}dA(w).$$
\end{remark}

\begin{theorem}
For any $\alpha>0$, let $d\mu_{\alpha}(w)=\frac{\alpha}{\pi}e^{-\alpha|w|^2}dA(w)$ be the weighted Gaussian measure. For any $\beta>2$ the operator $\mathcal{B}$ is bounded from $L^2(\mathbb{C},d\mu)$ into $L^2(\mathbb{C},d\mu_\beta)$. In particular, for any $f\in L^2(\mathbb{C},d\mu)$ it holds that
\begin{equation}
||\mathcal{B}(f)||_{L^2(\mathbb{C},d\mu_\beta)}\leq \frac{1}{\sqrt{\beta-2}} ||f||_{L^2(\mathbb{C},d\mu)}.
\end{equation}
\end{theorem}
\begin{proof}
For $z\in \mathbb{C}$ fixed, we have $$\mathcal{B}(f)(z)=e^{-|z|^2}\int_{\mathbb{C}} K(z,w)f(w)d\mu(w).$$
Thus, setting $K_z(w)=K(z,w)$ for any $w\in\mathbb{C}$ and noting that $f,K_z \in L^2(\mathbb{C},d\mu)$ we obtain
$$\displaystyle |\mathcal{B}(f)(z)|\leq e^{-|z|^2}\int_{\mathbb{C}}|K(z,w)||f(w)|d\mu(w),$$
and using the Cauchy Schwarz inequality
$$|\mathcal{B}(f)(z)|\leq e^{-|z|^2}||K_z||_{L^2(\mathbb{C},d\mu)}||f||_{L^2(\mathbb{C},d\mu)}.$$
Moreover, developing some calculations using Gaussian integrals we get

$$ ||K_z||_{L^2(\mathbb{C},d\mu)} =e^{2|z|^2},$$

and so
$$|\mathcal{B}(f)(z)|\leq e^{|z|^2} ||f||_{L^2(\mathbb{C},d\mu)}.$$
As a consequence, for any $\beta>2$ we have $$||\mathcal{B}(f)||^{2}_{L^2(\mathbb{C},d\mu_\beta)}\leq ||f||^{2}_{L^2(\mathbb{C},d\mu)}\left(\frac{1}{\pi}\int_{\mathbb{C}}e^{-(\beta-2)|z|^2}dA(z)\right)=\frac{1}{\beta-2}||f||^{2}_{L^2(\mathbb{C},d\mu)}.$$
Hence, $\mathcal{B}$ is a bounded operator from $L^2(\mathbb{C},d\mu)$ into $L^2(\mathbb{C},d\mu_\beta)$ with $\beta>2$.
\end{proof}
\begin{remark}
In particular, for any $\beta>2$ we have
$$||\mathcal{B}(f)||_{L^2(\mathbb{C},d\mu_\beta)}\leq ||\mathcal{B}(f)||_{\mathcal{SF}(\mathbb{C})}.$$
This explains somehow that $$\mathcal{SF}(\mathbb{C})\subset L^2(\mathbb{C},d\mu_\beta),\quad  \beta>2.$$
\end{remark}
Now, we can prove the following
\begin{lemma}\label{Naction}
For any $p,q=0, 1, 2, \ldots$, we have
$$ \mathcal{B}(H_{p,q})(z,\overline{z})=z^p\overline{z}^q.$$
We have also
   $$||\mathcal{B}(H_{p,q})||_{\mathcal{SF}(\mathbb{C})}=||H_{p,q}||_{L^2(\mathbb{C},d\mu)}.$$
\end{lemma}
\begin{proof}
For any $u,v\in\mathbb{C}$ we have (see \cite{Ismail, Ito} )
\begin{equation}
\displaystyle \sum_{m,n=0}^{\infty}H_{m,n}(z,\overline{z})\frac{u^mv^n}{m!n!}=e^{uz+v\overline{z}-uv}, \quad \forall z\in\mathbb{C}.
\end{equation}
For $w\in\mathbb{C}$ we set
\begin{equation}
\displaystyle \mathcal{R}(z,w):= \sum_{m,n=0}^{\infty}\frac{\overline{z}^mz^n}{m!n!} H_{m,n}(w,\overline{w})=e^{\overline{z}w+z\overline{w}-|z|^2}, \quad z\in\mathbb{C}.
\end{equation}
It is immediate that
$$\overline{\mathcal{R}(z,w)}=\mathcal{R}(z,w), \quad z,w\in\mathbb{C}.$$
From formula \eqref{Noperator} we obtain
 \[
    \begin{split}
   \displaystyle \mathcal{B}(H_{p,q})(z) &= \int_{\mathbb{C}} e^{-|z|^2} K(z,w)H_{p,q}(w,\overline{w})d\mu(w)\\
      &= \int_{\mathbb{C}} e^{-|z|^2+z\overline{w}+w\overline{z}}H_{p,q}(w,\overline{w})d\mu(w)\\
        &= \int_{\mathbb{C}}\left(\sum_{m,n=0}^{\infty}\frac{\overline{z}^mz^n}{m!n!} H_{m,n}(w,\overline{w}) \right) H_{p,q}(w,\overline{w})d\mu(w)
        \\
        &= \int_{\mathbb{C}}\mathcal{R}(z,w)H_{p,q}(w,\overline{w})d\mu(w)
        \\
        &= \int_{\mathbb{C}}\overline{\mathcal{R}(z,w)}H_{p,q}(w,\overline{w})d\mu(w)
         \\
        &= \sum_{m,n=0}^{\infty}\frac{z^m\overline{z}^n}{m!n!}\int_{\mathbb{C}}\overline{H_{m,n}(w,\overline{w})}H_{p,q}(w,\overline{w})d\mu(w)
            \\
        &= \sum_{m,n=0}^{\infty}\frac{z^m\overline{z}^n}{m!n!} \langle H_{p,q},H_{m,n} \rangle_{L^2(\mathbb{C},d\mu)}.
            \\
        &
          \end{split}
   \]

   However, we know that the complex Hermite polynomials $(H_{m,n})_{m,n\geq 0}$ form an orthonormal basis of $L^2(\mathbb{C},d\mu)$ (see \cite{Ito}) so that we have

   $$\langle H_{p,q},H_{m,n} \rangle_{L^2(\mathbb{C},d\mu)}=p!q!\delta_{(p,q);(m,n)}.$$

   In particular, this leads to
   $$\mathcal{B}(H_{p,q})(z,\overline{z})=z^p\overline{z}^q, \quad z\in\mathbb{C}.$$
   As a consequence, it is clear that for any $p,q= 0, 1, 2, \ldots$ we have

   $$||\mathcal{B}(H_{p,q})||_{\mathcal{SF}(\mathbb{C})}=||z^p\overline{z}^q||_{\mathcal{SF}(\mathbb{C})}=\sqrt{p!q!}=||H_{p,q}||_{L^2(\mathbb{C},d\mu)}.$$
\end{proof}
\begin{theorem}
The integral transform $\mathcal{B}$ is an isomorphism from $L^2(\mathbb{C},d\mu)$ onto $\mathcal{SF}(\mathbb{C})$ and for any $f\in L^2(\mathbb{C},d\mu)$ we have
\begin{equation}
||\mathcal{B}(f)||_{\mathcal{SF}(\mathbb{C})}=||f||_{L^2(\mathbb{C},d\mu)}.
\end{equation}
\end{theorem}
\begin{proof}
For any $f\in L^2(\mathbb{C},d\mu)$ we can write the following decomposition using Hermite polynomials
$$f(z)=\displaystyle \sum_{p,q=0}^{\infty}H_{p,q}(z,\overline{z})\alpha_{p,q}, \quad z\in\mathbb{C},$$
and $$||f||_{L^2(\mathbb{C},d\mu)}^{2}=\displaystyle \sum_{p,q=0}^{\infty}p!q!|\alpha_{p,q}|^2 .$$
Thus, thanks to Lemma \ref{Naction} we get
 \[
    \begin{split}
   \displaystyle \mathcal{B}(f)(z) &=  \sum_{p,q=0}^{\infty}\mathcal{B}(H_{p,q})(z,\overline{z})\alpha_{p,q} \\
      &=\sum_{p,q=0}^{\infty}z^p\overline{z}^q\alpha_{p,q}.
        \\
        &
          \end{split}
   \]

   Then, using also Proposition \ref{seqNF} we obtain
 \[
    \begin{split}
   \displaystyle || \mathcal{B}(f)||_{\mathcal{SF}(\mathbb{C})}^2&=  \sum_{p,q=0}^{\infty}p!q!|\alpha_{p,q}|^2 \\
      &=||f||^2_{L^2(\mu)}.\\
        &
        \\
        &
          \end{split}
   \]

   We observe that the surjectivity of the transform $\mathcal{B}$ is a direct consequence of Lemma \ref{Naction}. This allows to consider the integral transform $\mathcal{B}:L^2(\mathbb{C},d \mu)\longrightarrow \mathcal{SF}(\mathbb{C})$ which defines an unitary operator.

\end{proof}
\begin{proposition} \label{d^1}
For any $f\in L^2(\mathbb{C},d\mu)$, it holds that $$\displaystyle \frac{\partial}{\partial z}\mathcal{B}(f)(z)=-\overline{z}\mathcal{B}(f)(z)+\mathcal{B}(\overline{w}f)(z).$$
In a similar way, we have also
$$\displaystyle \frac{\partial}{\partial \overline{z}}\mathcal{B}(f)(z)=-z\mathcal{B}(f)(z)+\mathcal{B}(wf)(z).$$

\end{proposition}
\begin{proof}
We observe that thanks to the definition of $\mathcal{B}$ we have
$$\displaystyle \mathcal{B}(f)(z)=\frac{1}{\pi}\int_{\mathbb{C}} e^{-|z|^2}K(z,w)f(w)e^{-|w|^2}dA(w), \quad z\in\mathbb{C}.$$

Thus, we have
\begin{equation}\label{dzn}
\displaystyle \frac{\partial}{\partial z}\mathcal{B}(f)(z)=\frac{1}{\pi}\int_{\mathbb{C}}\frac{\partial}{\partial z} (e^{-|z|^2}K(z,w))f(w)e^{-|w|^2}dA(w), \quad z\in\mathbb{C}.
\end{equation}

Applying the Leibniz rule and Proposition \ref{Fockpptd} we get

 \[
    \begin{split}
   \displaystyle \frac{\partial}{\partial z} (e^{-|z|^2}K(z,w)) &= e^{-|z|^2}\frac{\partial}{\partial z} K(z,w)+\frac{\partial}{\partial z}(e^{-|z|^2})K(z,w) \\
      &=e^{-|z|^2}\overline{w}K(z,w)-\overline{z}e^{-|z|^2}K(z,w)\\
        &=(\overline{w}-\overline{z})e^{-|z|^2}K(z,w).
        \\
        &
          \end{split}
   \]

Hence, we insert the previous computations in the formula \eqref{dzn} and  this leads to 

 \[
    \begin{split}
   \displaystyle \frac{\partial}{\partial z} \mathcal{B}(f)(z)&= \frac{1}{\pi}\int_{\mathbb{C}} e^{-|z|^2}(\overline{w}-\overline{z})K(z,w)f(w)e^{-|w|^2}dA(w)\\
      &=\mathcal{B}(\overline{w}f)(z)-\overline{z}\mathcal{B}(f)(z).\\
        &
          \end{split}
   \]

\end{proof}
As a consequence of the previous result we prove
\begin{proposition}\label{Ndzbarz}
It holds that
\begin{equation}
\mathcal{B}^{-1}\left(\frac{\partial}{\partial z}+\overline{z}\right)\mathcal{B}=M_{\overline{w}},
\end{equation}
and \begin{equation}
\mathcal{B}^{-1}\left(\frac{\partial}{\partial \overline{z}}+z\right)\mathcal{B}=M_{w}.
\end{equation}
\end{proposition}
\begin{proof}
We note that thanks to Proposition \ref{d^1} we can write
 $$\displaystyle \frac{\partial}{\partial z}\mathcal{B}(f)(z)+\overline{z}\mathcal{B}(f)(z)=\mathcal{B}(\overline{w}f)(z).$$
Then, using the fact that $\mathcal{B}$ is an isomorphism it is easy to see that
$$\mathcal{B}^{-1}\left(\frac{\partial}{\partial z}+\overline{z}\right)\mathcal{B}=M_{\overline{w}}.$$
In a similar way we can prove the second statement.

\end{proof}
As a consequence, we can prove the following
\begin{corollary}
It holds that
\begin{equation}
\mathcal{B}^{-1}\left(\frac{1}{4}\Delta_z+\frac{\partial}{\partial\overline{z}} \overline{z}+z\frac{\partial}{\partial z}+|z|^2\right)\mathcal{B}=M_{|w|^2}.
\end{equation}
\end{corollary}
\begin{proof}
We observe that using the two expressions proved in Proposition \ref{Ndzbarz} we obtain that

$$\mathcal{B}^{-1}\left(\frac{\partial}{\partial z}+\overline{z}\right)\left(\frac{\partial}{\partial \overline{z}}+z\right)\mathcal{B}=M_{\overline{w}}M_w=M_{|w|^2}.$$
On the other hand, we have
$$\left(\frac{\partial}{\partial z}+\overline{z}\right)\left(\frac{\partial}{\partial \overline{z}}+z\right)=\left(\frac{1}{4}\Delta_z+\frac{\partial}{\partial\overline{z}} \overline{z}+z\frac{\partial}{\partial z}+|z|^2\right),$$
and this ends the proof.
\end{proof}
In the next result, we will calculate $\displaystyle \frac{\partial^n}{\partial \overline{z}^n}\mathcal{B}(f)(z)$ and $\displaystyle \frac{\partial^n}{\partial z^n}\mathcal{B}(f)(z)$:
\begin{proposition}
For any $n= 0, 1,2, \ldots$, it holds that
\begin{equation}\label{d^n}
\displaystyle \frac{\partial^n}{\partial z^n}\mathcal{B}(f)(z)=\sum_{k=0}^n(-1)^k {n \choose k}\overline{z}^k\mathcal{B}(\overline{w}^{n-k}f)(z)=\mathcal{B}((\overline{w}-\overline{z})^nf)(z),
\end{equation}

and

\begin{equation}
\displaystyle \frac{\partial^n}{\partial \overline{z}^n}\mathcal{B}(f)(z)=\sum_{k=0}^n(-1)^k {n \choose k}z^k\mathcal{B}(w^{n-k}f)(z)=\mathcal{B}((w-z)^nf)(z).
\end{equation}
\end{proposition}
\begin{proof}
The proof is done by induction. For $n=0$ the relations are true. We assume that \eqref{d^n} holds for some $n$ and we prove that we have

\begin{equation}\label{d^n+1}
\displaystyle \frac{\partial^{n+1}}{\partial z^{n+1}}\mathcal{B}(f)(z)=\sum_{k=0}^{n+1}(-1)^k {n+1 \choose k}\overline{z}^k\mathcal{B}(\overline{w}^{n+1-k}f)(z).
\end{equation}
Indeed, using Proposition \ref{d^1} we have
$$\displaystyle \frac{\partial}{\partial z}\mathcal{B}(f)(z)=-\overline{z}\mathcal{B}(f)(z)+\mathcal{B}(\overline{w}f)(z).$$
Thus, applying the operator $\displaystyle \frac{\partial^n}{\partial z^n}$ we get
$$\displaystyle \frac{\partial^{n+1}}{\partial z^{n+1}}\mathcal{B}(f)(z)=-\overline{z}\frac{\partial^n}{\partial z^n}\mathcal{B}(f)(z)+\frac{\partial^n}{\partial z^n}\mathcal{B}(\overline{w}f)(z).$$
Applying the induction hypothesis to the functions $f$ and $g=\overline{w}f$ we obtain the chain of equalities

\[
    \begin{split}
   \displaystyle \frac{\partial^{n+1}}{\partial z^{n+1}}\mathcal{B}(f)(z)
      &= -\overline{z}\left(\sum_{k=0}^{n}(-1)^k{n \choose k} \overline{z}^k\mathcal{B}(\overline{w}^{n-k}f)(z)\right)+\sum_{k=0}^{n}(-1)^k{n \choose k} \overline{z}^k\mathcal{B}(\overline{w}^{n+1-k}f)(z)\\
        &= \sum_{k=0}^{n}(-1)^{k+1}{n \choose k} \overline{z}^{k+1}\mathcal{B}(\overline{w}^{n-k}f)(z)+\sum_{k=0}^{n}(-1)^k{n \choose k} \overline{z}^k\mathcal{B}(\overline{w}^{n+1-k}f)(z)
        \\
        &= \sum_{h=1}^{n+1}(-1)^{h}{n \choose h-1} \overline{z}^{h}\mathcal{B}(\overline{w}^{n+1-h}f)(z)+\sum_{h=0}^{n}(-1)^h {n \choose h} \overline{z}^h\mathcal{B}(\overline{w}^{n+1-h}f)(z)
        \\
        &= {n \choose 0}\mathcal{B}(\overline{w}^{n+1}f)(z)+ \sum_{h=1}^{n}(-1)^{h}\left({n \choose h}+ {n \choose h-1} \right) \overline{z}^{h}\mathcal{B}(\overline{w}^{n+1-h}f)(z)
         \\
        &+(-1)^{n+1} {n \choose n}\overline{z}^{n+1}\mathcal{B}(f)(z).
          \end{split}
   \]
 By Pascal identity we have
   $${n \choose h}+ {n \choose h-1}={n+1 \choose h}, \quad \text{ for any } h\geq 1,$$
 so that
 $$\displaystyle \frac{\partial^{n+1}}{\partial z^{n+1}}\mathcal{B}(f)(z)=\sum_{k=0}^{n+1}(-1)^k {n+1 \choose k}\overline{z}^k\mathcal{B}(\overline{w}^{n+1-k}f)(z),$$
 and this ends the proof. The second statement follows using similar arguments.
\end{proof}

\begin{theorem}\label{Ndwbar'}
For any $f\in L^2(\mathbb{C},d\mu)$, it holds that  \begin{equation}
\mathcal{B}(\partial_{\overline{w}}f)(z)=\mathcal{B}(wf)(z)-z\mathcal{B}(f)(z).
\end{equation}
In a similar way, we have
\begin{equation}
\mathcal{B}(\partial_{w}f)(z)=\mathcal{B}(\overline{w}f)(z)-\overline{z}\mathcal{B}(f)(z).
\end{equation}
\end{theorem}
\begin{proof}
This result can be proved using some computations that are based on Proposition 7.2 of \cite{Shige}. Indeed, setting $d\mu(z):=\frac{1}{\pi}e^{-|z|^2}dA(z)$ we get
$$\displaystyle \int_{\mathbb{C}}u(w) \overline{(\partial_{\overline{w}}v(w))}d\mu(w)=\int_{\mathbb{C}}\left(-\partial_{w}+\overline{w}\right)(u(w))\overline{v(w)}d\mu(w).$$
In particular, this means that on $L^2(\mathbb{C},d\mu)$ we have

$$(\partial_w)^{*}=-\partial_w+\overline{w}.$$
Then, using this fact we deduce the following computations
\[
    \begin{split}
   \displaystyle M_z\mathcal{B}(f)(z) &=z\mathcal{B}(f)(z) \\
      &=e^{-|z|^2}\int_{\mathbb{C}} zK(z,w)f(w)d\mu(w)\\
        &=e^{-|z|^2}\int_{\mathbb{C}} \frac{\partial}{\partial \overline{w}}(K(z,w))f(w)d\mu(w)
        \\
        &=e^{-|z|^2}\langle \partial_{\overline{w}}K_z,\overline{f} \rangle_{\mu}
        \\
        &= e^{-|z|^2}\langle K_z,\partial_{\overline{w}}^{*}(\overline{f}) \rangle_{\mu}
         \\
        &= e^{-|z|^2}\langle K_z,(-\partial_{w}+\overline{w})(\overline{f}) \rangle_{\mu}
        \\
        &=e^{-|z|^2} \int_{\mathbb{C}} K(z,w)\overline{(-\partial_{w}+\overline{w})(\overline{f})(w)}d\mu(w)
        \\
        &=  e^{-|z|^2}\int_{\mathbb{C}} K(z,w) (-\partial_{\overline{w}}f+wf)d\mu(w)
        \\
        &=\mathcal{B}(-\partial_{\overline{w}}f+wf)(z)
        \\
        &=\mathcal{B}(wf)(z)-\mathcal{B}(\partial_{\overline{w}}f)(z),
          \end{split}
   \]

where we used
    $\overline{\partial_{w}f}=\partial_{\overline{w}}\overline{f},$ which leads to
   $$\overline{(-\partial_w+\overline{w})(\overline{f})}=(-\partial_{\overline{w}}f+wf).$$
   This ends the proof.
\end{proof}
As a consequence of the previous result, we prove that the operators $w-\partial_{\overline{w}}$, $\overline{w}-\partial_w$ are similar to $M_z$ and $M_{\overline{z}}$ respectively:
\begin{proposition}
It holds that
\begin{equation}
\mathcal{B}\left(w-\partial_{\overline{w}}\right)\mathcal{B}^{-1}=M_z,
\end{equation}
and

\begin{equation}
\mathcal{B}\left(\overline{w}-\partial_{w}\right)\mathcal{B}^{-1}=M_{\overline{z}}.
\end{equation}

\end{proposition}
\begin{proof}
By Theorem \ref{Ndwbar'} we have
$$\mathcal{B}(\partial_{\overline{w}}f)(z)=\mathcal{B}(wf)(z)-z\mathcal{B}(f)(z),$$
from which we deduce
$$\mathcal{B}((w-\partial_{\overline{w}})f)(z)=z\mathcal{B}(f)(z).$$
Hence, using the fact that the integral transform $\mathcal{B}$ is a unitary operator from $L^2(\mathbb{C},d\mu)$ onto $\mathcal{SF}(\mathbb{C})$ we multiply the operator $\mathcal{B}^{-1}$ on the right and obtain $$\mathcal{B}\left(w-\partial_{\overline{w}}\right)\mathcal{B}^{-1}=M_z.$$

The second part of the statement can be proved similarly.
\end{proof}
\begin{proposition}\label{p5.12}
We have
\begin{equation}
\partial_z \mathcal{B}(f)=\mathcal{B}(\partial_w f)
\end{equation}
and
\begin{equation}
\partial_{\overline{z}}N(f)=\mathcal{B}(\partial_{\overline{w}}f).
\end{equation}
\end{proposition}
\begin{proof}
Applying Proposition \ref{d^1} and Theorem \ref{Ndwbar'}  we have
$$\displaystyle \partial_z\mathcal{B}(f)(z)=-\overline{z}\mathcal{B}(f)(z)+\mathcal{B}(\overline{w}f)(z)$$
and $$\mathcal{B}(\partial_{w}f)(z)=\mathcal{B}(\overline{w}f)(z)-\overline{z}\mathcal{B}(f)(z).$$
Therefore, from the two previous relations we obtain

$$\displaystyle \partial_z\mathcal{B}(f)(z)=\mathcal{B}(\partial_{w}f)(z).$$
The second statement can be proved in a similar way and this will end the proof.
\end{proof}

In the table below, we summarize different operators that are equivalent to each others thanks to the Berezin transform $\mathcal{B}$.
\begin{table}[ht]
\caption{Equivalent operators using the transform $\mathcal{B}$} 
\centering 
\large
\begin{tabular}{cc } 
\hline 
 $\mathcal{SF}(\mathbb{C})$ & $L^2(\mathbb{C},d\mu)$  \\ [1ex] 

\hline 
\vspace{0.2cm}
$\frac{\partial}{\partial z}+\overline{z}$ & $\overline{w}$  \\  
\vspace{0.2cm}
$\frac{\partial}{\partial\overline{z}}+z$ & $w$  \\
\vspace{0.2cm}
$\frac{\partial^n}{\partial z^n}$ & $\frac{\partial^n}{\partial w^n}$\\
\vspace{0.2cm}
$\frac{\partial^n}{\partial\overline{z}^n}$ & $\frac{\partial^n}{\partial \overline{w}^n}$ \\
\vspace{0.2cm}
$\frac{1}{4}\Delta_z+\frac{\partial}{\partial\overline{z}} \overline{z}+z\frac{\partial}{\partial z}+|z|^2$ & $|w|^2$  \\
\vspace{0.2cm}
$M_{|z|^2}$ & $\frac{1}{4}\Delta_w-\frac{\partial}{\partial\overline{w}} \overline{w}-w\frac{\partial}{\partial w}+|w|^2$   \\
\vspace{0.2cm}
$z$ & $w-\frac{\partial}{\partial\overline{w}}$  \\
\vspace{0.2cm}
$\overline{z}$ & $\overline{w}-\frac{\partial}{\partial w}$  \\
\hline 
\end{tabular}
\label{table:nonlin} 
\end{table} \newpage
Inspired by \cite{Zhu2} we consider the following integral operator on the Fock space $\mathcal{F}(\mathbb{C})$
\begin{definition}
Let $\varphi:\mathbb{C}\longrightarrow \mathbb{C}$ and denote by $d\mu(w)=\frac{1}{\pi}e^{|w|^2}dA(w)$ the normalized Gaussian measure. Then, we define the following integral transform when it exists
\begin{equation}
\mathcal{S}_\varphi(f)(z):=\displaystyle \int_{\mathbb{C}}e^{z\overline{w}}f(w)\varphi(w-z)d\mu(w),\quad f\in\mathcal{F}(\mathbb{C}).
\end{equation}
\end{definition}
\begin{remark}
 If $\varphi=1$, it is clear by the reproducing kernel property for the Fock space that in this case $$ \mathcal{S}_\varphi(f)(z):=\displaystyle \int_{\mathbb{C}}e^{z\overline{w}}f(w)d\mu(w)=f(z),\quad \forall f\in\mathcal{F}(\mathbb{C}).$$
\end{remark}
\begin{example}
For every $a\in\mathbb{C}$ set $\varphi_a(w)=e^{\overline{a}w}$. Then, we have
$$\mathcal{S}_{\varphi_z}(f)(z)=\displaystyle \int_{\mathbb{C}}e^{z\overline{w}}f(w)\varphi_z(w-z)d\mu(w)=\mathcal{B}(f)(z).$$
 It turns out that the transform $\mathcal{B}$ is a particular case of the general integral operator $\mathcal{S}_{\varphi}$.
\end{example}

\section{The polyanalytic Hardy space of infinite order and Gleason problem}
\setcounter{equation}{0}
Let us denote by $\mathbb{D}$ the unit disk and by $\mathcal{H}^2(\mathbb{D})$ the classical Hardy space. In this section, we will prove the following main result:
\begin{theorem}
The reproducing kernel Hilbert space with reproducing kernel $\displaystyle \frac{1}{(1-z\overline{w})(1-\z w)}$ is, up to a multiplicative positive factor, the only reproducing kernel
Hilbert space of polyanalytic functions of infinite order, regular at the origin, and for which
\begin{eqnarray}
R_\infty^*&=&M_z\\
L_\infty^*&=&M_{\overline{z}}.
 \end{eqnarray}
\end{theorem}


To prove this theorem we need a couple of preliminary results, including a sequential characterization of the space $\mathcal{H}(\mathsf{K})$.
\begin{definition}
The polyanalytic Hardy space of infinite order $\mathcal{SH}(\mathbb{D})$
 is the space of functions of the form \begin{equation}
f(z)=\displaystyle \sum_{n=0}^{\infty}\overline{z}^nf_n(z),
\end{equation}
satisfying
\begin{enumerate}
\item[i)] $f_n\in \mathcal{H}^2(\mathbb{D})$ for any $n\geq 0$;
\item[ii)]$\displaystyle ||f||^{2}_{\mathcal{SH}(\mathbb{D})}= \sum_{n=0}^{\infty} ||f_n||^{2}_{\mathcal{H}^2(\mathbb{D})}<\infty.$
\end{enumerate}
Then, we consider the scalar product on $\mathcal{SH}(\mathbb{D})$ given by \begin{equation}\displaystyle
\langle f,g \rangle_{\mathcal{SH}(\mathbb{D})}:=\sum_{k=0}^{\infty}\langle f_k,g_k \rangle_{\mathcal{H}^2(\mathbb{D})},
\end{equation}
 for any $f=\displaystyle \sum_{k=0}^{\infty}\overline{z}^kf_k$ and $g=\displaystyle \sum_{k=0}^{\infty}\overline{z}^kg_k$ with $f_k,g_k\in\mathcal{H}^2(\mathbb{D})$ for every $k\geq 0$.
\end{definition}

\begin{proposition}
A function $f:\mathbb{D}\longrightarrow \mathbb{C}$ belongs to $\mathcal{SH}(\mathbb{D})$ if and only if $f$ is of the form $$\displaystyle f(z)=\sum_{(m,n)\in\mathbb{N}^2}z^m\overline{z}^n\alpha_{m,n},$$
with $(\alpha_{m,n})\subset\mathbb{C}$ and such that \begin{equation}\displaystyle
||f||^{2}_{\mathcal{SH}(\mathbb{D})}=\sum_{(m,n)\in\mathbb{N}^2}|\alpha_{m,n}|^2<\infty.
\end{equation}
Moreover, if for any $(m,n)\in\mathbb{N}^2$ we set $\upsilon_{m,n}(z,\overline{z})=\displaystyle z^m\bar{z}^n$. Then, the family of functions $\lbrace \upsilon_{m,n} \rbrace_{m,n\geq 0}$ form an orthonormal basis of $\mathcal{SH}(\mathbb{D})$.
\end{proposition}
\begin{proof}
This proof follows the arguments used to prove Proposition \ref{seqNF}.
\end{proof}

\begin{lemma}
We have
$$\mathcal{H}(\mathsf{K})=\mathcal{SH}(\mathbb{D}).$$
Moreover, it holds that

\begin{equation}
\mathsf{K}(z,w)=\displaystyle \sum_{m,n=0}^{\infty}\upsilon_{m,n}(z,\overline{z})\overline{\upsilon_{m,n}(w,\overline{w})},
\end{equation}
for every $z,w\in\mathbb{D}$.
\end{lemma}
\begin{proof}
We note that $(\upsilon_{m,n})_{m,n\geq 0}$ is an orthonormal basis of the space $\mathcal{SH}(\mathbb{D})$. Thus, the associated reproducing kernel is given by the following series which converges uniformly on each compact so that
$$\displaystyle \sum_{m,n=0}^{\infty}\upsilon_{m,n}(z,\bar{z})\overline{\upsilon_{m,n}(w,\bar{w})}<\infty, \text{ for any } z,w\in\mathbb{D}.$$
More precisely, for any $(z,w)\in\mathbb{D}^2$ we have the equalities
\[
    \begin{split}
   \displaystyle  \sum_{m,n=0}^{\infty}\upsilon_{m,n}(z,\bar{z})\overline{\upsilon_{m,n}(w,\bar{w})}
      &= \sum_{m,n=0}^{\infty} z^m\bar{z}^n \bar{w}^m w^n \\
        &= \left(\sum_{m=0}^{\infty} z^m\bar{w}^m \right) \left(\sum_{n=0}^{\infty} w^n\bar{z}^n\right)
        \\
        &= \frac{1}{(1-z\overline{w})} \frac{1}{(1-\z w)}
        \\
        &=\frac{1}{(1-z\overline{w})(1-\z w)}
         \\
        &=\mathsf{K}(z,w).
          \end{split}
   \]
\end{proof}

\begin{lemma}\label{HardyR_0}
It holds that
\begin{equation}\displaystyle
\langle R_\infty(f),g\rangle_{\mathcal{SH}(\mathbb{D})}=\langle f,M_{z}g\rangle_{\mathcal{SH}(\mathbb{D})},
\end{equation}
and \begin{equation}\displaystyle
\langle L_\infty(f),g\rangle_{\mathcal{SH}(\mathbb{D})}=\langle f,M_{\overline{z}}g\rangle_{\mathcal{SH}(\mathbb{D})}.
\end{equation}
\end{lemma}
\begin{proof}
Let $\displaystyle f=\sum_{k=0}^{\infty}\overline{z}^kf_k$ and $\displaystyle g=\sum_{k=0}^{\infty}\overline{z}^kg_k$ in $\mathcal{SH}(\mathbb{D})$. Since we have $$\displaystyle M_z(g)=\sum_{k=0}^{\infty}\overline{z}^kM_z(g_k),$$ it follows that
\[
    \begin{split}
   \displaystyle  \langle R_\infty(f),g\rangle_{\mathcal{SH}(\mathbb{D})}
      &= \sum_{k=0}^{\infty}\langle R_0(f_k),g_k \rangle_{\mathcal{H}^2(\mathbb{D})}\\
        &= \sum_{k=0}^{\infty}\langle f_k,(R_0)^*g_k \rangle_{\mathcal{H}^2(\mathbb{D})}
        \\
        &= \sum_{k=0}^{\infty}\langle f_k,M_z(g_k) \rangle_{\mathcal{H}^2(\mathbb{D})}
        \\
        &= \sum_{k=0}^{\infty}\langle f_k,zg_k \rangle_{\mathcal{H}^2(\mathbb{D})}
         \\
        &=\langle f,M_{z}(g)\rangle_{\mathcal{SH}(\mathbb{D})}.
          \end{split}
   \]
  In a similar way, we can prove the second part of the statement.
\end{proof}
   \begin{example}
  On  the bidisc $\mathbb{D}^2$ we consider the inner function defined by $$j(z_1,z_2)=\displaystyle \frac{z_1+z_2+2z_1z_2}{z_1+z_2+2}, \quad \forall (z_1,z_2)\in\mathbb{D}^2.$$
  Then, $j(z_1,z_2)$ is a contractive multiplier of the Hardy space $\mathcal{H}^2(\mathbb{D}^2)$ and hence (dimension 2) is in the Schur-Agler class, see \cite{BB2012,BVD2015}. Moreover, if we set
$$\rho_w(z)=1-z\overline{w},$$
we can consider the kernel function
  \begin{equation}
  \mathsf{K}_j((z_1,z_2);(w_1,w_2)):= \frac{1-j(z_1,z_2)\overline{j(w_1,w_2)}}{\rho_{w_1}(z_1)\rho_{w_2}(z_2)}, \quad \forall (z_1,z_2); (w_1,w_2) \in\mathbb{D}^2.
  \end{equation}
  We note that

\[
    \begin{split}
   \displaystyle \mathsf{K}_j((z_1,z_2);(w_1,w_2))
      &=\frac{2}{(z_1+z_2+2)(\overline{w_1}+\overline{w_2}+2)}\\
        & \cdot \left( \frac{(z_1+1)(\overline{w_1}+1)}{(1-z_1\overline{w_1})}+\frac{(z_2+1)(\overline{w_2}+1)}{(1-z_2\overline{w_2})} \right).
        \\
        &
          \end{split}
   \]

Then, by taking $z_1=z, z_2=\overline{z}$ and $w_1=w, w_2=\overline{w}$ we have
\begin{equation}\label{j}
\displaystyle j(z,\overline{z})=\frac{z+\overline{z}+2|z|^2}{z+\overline{z}+2}=\frac{\Re(z)+|z|^2}{1+\Re(z)}, \quad \forall z\in\mathbb{D},
\end{equation}

and we can write
$$j(z,\overline{z})=\frac{P(z,\overline{z})}{Q(z,\overline{z})},$$
where both the polynomials $P$ and $Q$ are polyanalytic of order $2$. We observe that $Q(z,\overline{z})=0$ if and only if $\Re(z)=-1,$ so, $Q(z,\overline{z})\neq 0$ for every $z\in \mathbb{D}$.
On the other hand, we note that $j(z,\overline{z})=1$ on the boundary $\partial\mathbb{D}$.
We have

\begin{equation}
\displaystyle \mathsf{K}_j(z,w)=\frac{1-j(z,\overline{z})\overline{j(w,\overline{w})}}{\rho_w(z)\rho_{\overline{w}}(\overline{z})}, \quad \forall (z,w) \in\mathbb{D}^2,
\end{equation}
and, as a consequence,

\[
    \begin{split}
   \displaystyle \mathsf{K}_j((z,\overline{z});(w,\overline{w}))
      &=\frac{2}{(z+\overline{z}+2)(\overline{w}+w+2)}\\
        & \cdot \left( \frac{(z+1)(\overline{w}+1)}{(1-z\overline{w})}+\frac{(\overline{z}+1)(w+1)}{(1-\overline{z}w)} \right).
        \\
        &
          \end{split}
   \]
   It is important to note also that the function $j$ given by \eqref{j} is polyrational in the sense of \cite[pp 175]{Balk1991}.
    \end{example}
    \begin{remark}
According to \cite{agler2, BT1998} we observe that \begin{equation}
    j(z_1,z_2)=C(I_2-ZA)^{-1}ZB,
    \end{equation}
   and using the unitary matrix $M$ given by

\[
  M=\begin{pmatrix}A & B \\
    C & D
 \end{pmatrix}=\begin{pmatrix}-\frac{1}{2}&-\frac{1}{2}&\frac{1}{\sqrt{2}}\\
    -\frac{1}{2}&-\frac{1}{2}&-\frac{1}{\sqrt{2}}\\
  \frac{1}{\sqrt{2}}&-\frac{1}{\sqrt{2}}&0
 \end{pmatrix},
\]
we have
\[
  A=\begin{pmatrix}-\frac{1}{2}&-\frac{1}{2}\\
    -\frac{1}{2}&-\frac{1}{2}
 \end{pmatrix},  C=\begin{pmatrix} \frac{1}{\sqrt{2}}&-\frac{1}{\sqrt{2}}
 \end{pmatrix}, B=\begin{pmatrix}\frac{1}{\sqrt{2}}\\
    -\frac{1}{\sqrt{2}}
 \end{pmatrix}
\]
and $D=0$ where we set
$
Z=\begin{pmatrix}z_1& 0 \\
    0& z_2
 \end{pmatrix}.
$
    \end{remark}

Inspired by \cite{ABDS2001} we can prove the following result concerning the backward shift operators considered in Definition \ref{RinfinityL}  in the case of the polyanalytic Hardy space of infinite order $\mathcal{SH}(\mathbb{D})$. Indeed, we have
\begin{theorem}
A function $f\in\mathcal{SH}(\mathbb{D})$ is a common eigenfunction for the backward shift operators $R_\infty$ and $L_\infty$ with corresponding eigenvalues given respectively by $\lambda_1, \lambda_2\in\mathbb{D}$ if and only if

\begin{equation}\label{fexpression}
\displaystyle
f(z,\overline{z})=\frac{f(0,0)}{(1-\lambda_1 z)(1-\lambda_2\overline{z})}, \quad \forall z\in\mathbb{D}.
\end{equation}

\end{theorem}
 \begin{proof}
It is easy to check that if $f$ is of the form \eqref{fexpression}, then we have
$R_\infty(f)=\lambda_1f$ and $L_\infty(f)=\lambda_2f.$
For the converse, we assume that $f$ is an eigenfunction for $R_\infty$ and $L_\infty$ with corresponding eigenvalues given by $\lambda_1$ and $\lambda_2$. By writing $$\displaystyle f(z,\overline{z})=\sum_{n=0}^\infty\overline{z}^nf_n(z),$$
we observe that $R_\infty(f)=\lambda_1 f$ if and only if $$\displaystyle \sum_{n=0}^{\infty} \overline{z}^nR_0(f_n)(z)=\sum_{n=0}^\infty \overline{z}^n(\lambda_1f_n)(z).$$

In particular, this shows that  $$R_0(f_n)(z)=\lambda_1 f_n(z), \quad \text{ for any } n= 0,1, 2, \ldots.$$
Using the classical result on the backward shift operator $R_0$ for any $n= 0, 1, 2, \ldots$ we get
$$\displaystyle f_n(z,\overline{z})=\frac{f_n(0)}{1-\lambda_1z},\quad z\in\mathbb{D}.$$
Now, we insert the expression of $f_n$ in the $f$ decomposition which leads to the following calculations
\[
    \begin{split}
   \displaystyle f(z)
      &= \sum_{n=0}^{\infty}\overline{z}^nf_n(z)\\
        &= \sum_{n=0}^{\infty}\overline{z}^n\frac{f_n(0)}{1-\lambda_1z}
        \\
        &= \left(\sum_{m=0}^{\infty}z^m \lambda_1^m \right)\left(\sum_{n=0}^{\infty}\overline{z}^nf_n(0)\right)
         \\
        &=\sum_{m=0}^{\infty}z^mg_m(\overline{z})
          \end{split}
   \]
   where we set $\displaystyle g_m(\overline{z})= \lambda_1^m\left(\sum_{n=0}^{\infty}\overline{z}^nf_n(0)\right)$ for any $z\in\mathbb{D}$. On the other hand, we note that by definition $$\displaystyle L_\infty(f)=\sum_{m=0}^{\infty}z^mL_0(g_m),$$
   with $\displaystyle L_0(g_m)(\overline{z})=\frac{g_m(\overline{z})-g_m(0)}{\overline{z}}$. Then, following a similar reasoning as we did in the case of $R_\infty$,  using the fact that $L_\infty(f)=\lambda_2 f$ we deduce tha $L_0(g_m)(\overline{z})=\lambda_2 g_m(\overline{z})$. Thus, in particular this allows to write
   $$g_m(\overline{z})=\frac{g_m(0)}{(1-\lambda_2\overline{z})}, \quad z\in\mathbb{D}.$$
   So, now we insert $g_m$ in the expression of $f$ and get
   $$\displaystyle f(z)=\sum_{m=0}^\infty z^m \frac{g_m(0)}{1-\lambda_2\overline{z}}.$$
  We note that $g_m(0)=\lambda_1^mf_0(0)=\lambda_1^mf(0)$. Then, we obtain

  $$f(z)=\frac{f(0)}{1-\lambda_2\overline{z}}\sum_{m=0}^{\infty}\lambda_1^mz^m, \quad z\in \mathbb{D}.$$
  Hence, we conclude that

  $$\displaystyle
f(z,\overline{z})=\frac{f(0,0)}{(1-\lambda_1 z)(1-\lambda_2\overline{z})}, \quad \text{ for any } z\in\mathbb{D}.$$
 \end{proof}

\begin{lemma}
For any $f\in \mathcal{C}^1$ we have
\[
\frac{\rm d}{\rm dt}f(tx,ty)=z\partial f+\z \overline{\partial}f
\]
\end{lemma}

\begin{proof}
It follows from
  \[
    \begin{split}
      \frac{\rm d}{\rm dt}f(tx,ty)&=x\frac{\partial f  }{\partial x}(tx,ty)+y\frac{\partial f  }{\partial y}(tx,ty)\\
      &=\frac{1}{2}(x+iy)\left(\frac{\partial f  }{\partial x}(tx,ty)--i\frac{\partial f  }{\partial y}(tx,ty)\right)+\frac{1}{2}(x-iy)\left(\frac{\partial f  }{\partial x}(tx,ty)
        +i\frac{\partial f  }{\partial y}(tx,ty)\right)\\
      &=
      z\partial f+\z \overline{\partial}f.
    \end{split}
    \]
\end{proof}
Let
\begin{equation}
  f(z,\z)=zf_1(z,\z)+\z f_2(z,\z)
\end{equation}
where $f_1$ and $f_2$ are required to be in the same space as $f$. Then
\begin{eqnarray}
  (A_0 f)(z,\z)&=&\int_0^1\frac{\partial}{\partial z}f(tz,t\z)dt\\
  (B_0 f)(z,\z)&=&\int_0^1\frac{\partial}{\partial \z}f(tz,t\z)dt
\end{eqnarray}
We have
\[
    \frac{\rm d}{{\rm d}t}f(tz,t\z)=z\frac{\partial}{\partial z}f(tz,t\z)dt\\+
  \z \frac{\partial}{\partial \z}f(tz,t\z)dt
\]
and so
\begin{equation}
  \label{345}
  f(z,\z)-f(0,0)=zA_0f(z,\z)+\z B_0f(z,\z)
\end{equation}

\begin{lemma}
Let $\mathfrak M$ be a finite dimensional space in which \eqref{345} holds. Any $f\in \mathfrak{M}$ which is a common eigenfunction of $A_0$ and $B_0$ can be written in the form
\[
  f(z,\z)=\frac{f(0,0)}{1-az-b\z}.
  \]
\end{lemma}

\begin{proof}
  Since $A_0$ and $B_0$ commute they can be simultaneously triangularized.
  Let
  \[
    A_0f=a f\quad{\rm and}\quad B_0f=b f
  \]
  We have
  \[
  f(z,\z)=f(0,0)+(az+b\z)f
\]
and so
\[
  f(z,\z)=\frac{f(0,0)}{1-az-b\z}
  \]
\end{proof}

\begin{lemma}
The following equalities hold:
  \begin{eqnarray}
    A_0(z^n\z^m)&=&\frac{n}{n+m}z^{n-1}\z^m\\
    B_0(z^n\z^m)&=&\frac{m}{m+n}z^n\z^{m-1}.
\end{eqnarray}
\end{lemma}

\begin{proof}
It follows from
  \[
  \begin{split}
    A_0(z^n\z^m)&=\int_0^1(nz^{n-1}\z^m)(tz,t\z)dt\\
    &=nz^{n-1}\z^m\int_0^1nt^{n+m-1}dt\\
    &=\frac{n}{n+m}z^{n-1}\z^m,
  \end{split}
  \]
  and similarly for $B_0$.
 \end{proof}
\begin{remark}
We remark that $A_0$ is not $R_0$, in general; it will reduce to $R_0$ when $f$ is analytic. We note that both the operators $R_\infty$ and $A_0$ extend the classical backward shift operator $R_0$ on the Hardy space, but these two operators are different.
\end{remark}
\section{The Drury-Arveson space case}
\setcounter{equation}{0}
Let us consider the kernel function given by
\begin{equation}
  \label{arve}
k(z,w)=  \frac{1}{1-(z\overline{w}+\z w)}.
\end{equation}
We denote by $\mathfrak H(k)$ the associated reproducing kernel Hilbert space. Setting $z=x+iy$ and $w=t+iu$, we have
\begin{equation}
  \frac{1}{1-(z\overline{w}+\z w)}=\frac{1}{1-2(xt+yu)}
\end{equation}

We note that \eqref{arve} is a complete Nevanlinna-Pick kernel, meaning that
\[
\frac{1}{k(z,w)}=1-2xt-2yu
\]
has one positive square in $B(0,1/\sqrt{2})$. Such kernels were introduced by Agler, see \cite{agler} and also the paper of Quiggin \cite{PQ}.
\begin{lemma}\label{kresult}
The function \eqref{arve} is positive definite in $|z|<1/\sqrt{2}$ and the functions
\[
\frac{\partial k}{\partial t}\quad and\quad \frac{\partial k}{\partial u}
\]
belong to $\mathfrak H(k)$. Furthermore, it holds that
\begin{eqnarray}
  \langle f, \frac{\partial f}{\partial t}\rangle_{\mathfrak H(k)}&=&\frac{\partial f}{\partial x}\\
  \langle f, \frac{\partial f}{\partial u}\rangle_{\mathfrak H(k)}&=&\frac{\partial f}{\partial y}.
                                                                      \end{eqnarray}
\end{lemma}

\begin{proof}

For a fixed choice of $(t,u)$ and for $h\in\mathbb R$ small enough we set
  \[
f_h(x,y,t,u)=\frac{k(x,y,t,u+h)-k(x,y,t,u)}{h}.
\]
Then $f_h\in\mathfrak H(k)$ and
\[
\|f_h\|^2_{\mathcal H(k)}=\frac{k(t,u+t,t,u+h)+k(t,u,t,u)-2k(t,u+h,t,u)}{h^2}
\]
uniformly bounded in $h$ for $h$ small. Thus $f_h$ has a weakly convergent subsequence, with limit say $g_{t,u}$.
Since weak convergence implies pointwise convergence we have
\[
\begin{split}
  g_{t,u}(x,y)&=\langle g_{t,u},k(\cdot,\cdot, x,y)\rangle_{\mathfrak H(k)}\\
  &=\lim_{h\rightarrow 0} \langle f_h,k(\cdot,\cdot, x,y)\rangle_{\mathfrak H(k)}\\
  &=\lim_{h\rightarrow 0} \left\langle   \frac{k(\cdot,\cdot,t,u+h)-k(\cdot,\cdot,t,u)}{h},k(\cdot,\cdot, x,y)\right\rangle_{\mathfrak H(k)}\\
  &=\lim_{h\rightarrow 0}  \frac{k(x,y,t,u+h)-k(x,y,t,u)}{h}\\
  &= \frac{\partial f}{\partial u}(x,y,t,u).
  \end{split}
\]


Furthermore, for $f\in\mathfrak H(k)$, we have:
\[
  \begin{split}
    \langle f,g_{t,u}\rangle_{\mathfrak H(k)}=&=\lim_{h\rightarrow 0} \langle f,f_h\rangle_{\mathfrak H(k)}\\
    &=\lim_{h\rightarrow 0} \left\langle f(\cdot,\cdot),  \frac{k(\cdot,\cdot,t,u+h)-k(\cdot,\cdot,t,u)}{h}\right\rangle_{\mathfrak H(k)}\\
    &=\lim_{h\rightarrow 0} \frac {f(t,u+h)-f(t,u)}{h}\\
  &=\frac{\partial f}{\partial y}(t,u)
  \end{split}
\]
The other claims are proved in the same way.
\end{proof}

Iterating the above result we get:

\begin{corollary}
Let $k$ be as in \eqref{arve}, then for $n,m=0,1,2,\ldots$
  \begin{equation}
\frac{\partial^{n+m}k}{\partial t^n\partial u^m}(\cdot,\cdot, t,u)\in\mathfrak H(k)
\end{equation}
and
\begin{equation}
  \langle f(\cdot,\cdot), \frac{\partial^{n+m}k}{\partial t^n\partial u^m}(\cdot,\cdot, t,u)\rangle_{\mathfrak H(k)}=
  \frac{\partial^{n+m}f}{\partial x^n\partial y^m}(t,u).
\end{equation}
\end{corollary}
\begin{proof}
This a direct consequence of Lemma \ref{kresult}.
\end{proof}

\begin{corollary}
For for $n,m=0,1,2,\ldots$ we have
  \[
    x^ny^m\in\mathfrak H(k).
    \]
  \end{corollary}

  \begin{proof}
It suffices to set $t=u=0$ in the previous corollary.
    \end{proof}

    We now give a characterization of the space $\mathfrak H(k)$.

    \begin{proposition}
    The space  $\mathfrak H(k)$ consists of the functions of the form
\begin{equation}
  f(z,\z)=\sum_{a,b=0}^\infty c_{a,b}z^a\z^b
\end{equation}
with norm
\begin{equation}
  \|f\|^2=\sum_{a,b=0}^\infty |c_{a,b}|^2\frac{a!b!}{(a+b)!}
  \end{equation}
      \end{proposition}
    \begin{proof}
    It suffices to observe that we have
\begin{equation}
 k(z,w)= \sum_{a,b\in\mathbb N_0}\frac{(a+b)!}{a!b!}z^a\z^bw^b\overline{w}^{a}.
\end{equation}
\end{proof}
\begin{proposition}
  The operators $M_z$ and $M_{\overline{z}}$ are bounded in $\mathfrak H(k)$, with $\|M_{z}\|\le 1$ and $\|M_{\overline{z}}\|\le 1$. Their adjoints are given by
  \begin{eqnarray}
    M_z^*&=&A_0\\
    M_{\overline{z}}^*&=&B_0.
\end{eqnarray}
\end{proposition}
\begin{proof}
  The first claim follows from
  \[
\frac{1-z\overline{w}}{1-2{\rm Re}\, z\overline{w}}=1+\frac{\overline{z}w}{1-2{\rm Re}\, z\overline{w}}\ge 0
  \]
and
  \[
\frac{1-\overline{z}w}{1-2{\rm Re}\, z\overline{w}}=1+\frac{z\overline{w}}{1-2{\rm Re}\, z\overline{w}}\ge 0
\]
The second claim follows from
  \[
    \begin{split}
\langle M_z(z^n\z^m),z^u\z^v\rangle&=\langle z^{n+1}\z^m,z^u\z^v\rangle
\end{split}
\]
and similarly for $M_{\overline{z}}$.
\end{proof}

Following \cite[Corollary 2.4, p. 7]{AT} we introduce
\begin{eqnarray}
  (A_af)(z)&=&\frac{z}{1-2{\rm Re}\, z\overline{a}}f(z)\\
  (B_af)(z)&=&\frac{\overline{z}}{1-2{\rm Re}\, \overline{z}a}f(z)
               \end{eqnarray}
               with $a\in B(0,1/\sqrt{2})$.
               \begin{proposition}
                 Let $a\in B(0,1/\sqrt{2})$. The operators $A_a$ and $B_a$ are bounded and it holds that
                 \begin{equation}
                   \label{wer}
f(z)-f(a)=(z-a)(A_a^*f)(z)+(\overline{z}-\overline{a})(B_a^*f)(z),\quad f\in\mathfrak H(k).
                   \end{equation}
               \end{proposition}

               \begin{proof}
We have
                 \begin{eqnarray}
                   (A_a^*k_b)(z)&=\frac{\overline{b}}{1-2{\rm Re}\, (b\overline{a})}k_b(z)\\
                   (B_a^*k_b)(z)&=\frac{b}{1-2{\rm Re}\, (b\overline{a})}k_b(z)
                 \end{eqnarray}

                 \[
                   \begin{split}
                     (z-a)(A_a^*k_b)(z)+(\overline{z}-\overline{a})(B_a^*k_b)(z)&=\frac{(z-a)\overline{b}+(\overline{z}-\overline{a})b}{(1-2{\rm Re}\, b\overline{a})(1-2{\rm Re}\, z\overline{b})}\\
                     &=\frac{z\overline{b}+\overline{z}b-a\overline{b}-\overline{a}b}{(1-2{\rm Re}\, b\overline{a})(1-2{\rm Re}\, z\overline{b})}\\
                     &=k_b(z)-k_b(a)
                   \end{split}
                   \]
     \end{proof}

\begin{example}
  Let us consider the coefficients $c_n$ such that

\begin{equation}\label{cn}
 \displaystyle 1-\sqrt{1-t}=\sum_{n=1}^{\infty} c_nt^n, \quad t<1.
\end{equation}

Then,  the function
  \[
f_m(z,\z)=z+\sum_{n=1}^mc_n\z^{2m}
\]
is bounded by one in modulus in $|z|<\frac{1}{\sqrt{2}}$, but is not a Schur multiplier.
\end{example}
Indeed, recalling that the $c_n>0$ and satisfy $\sum_{n=1}^\infty c_n=1$ we have
\[
  \begin{split}
    \left|z+\sum_{n=1}^mc_n\z^{2n}\right|&\le |z|+\sum_{n=1}^m c_n|z|^{2n}\\
    &\le |z|+\sum_{n=1}^\infty c_n|z|^{2n}\\
    &=|z|+1-\sqrt{1-|z|^2}\\
    &\le \frac{1}{\sqrt{2}}+1-\sqrt{\frac{1}{2}}\\
      &=1.
  \end{split}
\]
But $\|f_m\|^2=1+\sum_{n=1}^mc_n\|\z^{2n}\|^2>1$.

For $a\in B(0,1/\sqrt{2})$ we set
\begin{equation}
  b_a(z)=\frac{(1-2|a|^2)\begin{pmatrix}z-a&\overline{z}-\overline{a}\end{pmatrix}}{1-{\rm Re}\, z\overline{a}}\sqrt{I_2-\begin{pmatrix}\overline{a}\\ a\end{pmatrix}
        \begin{pmatrix}a&\overline{a}\end{pmatrix}}
\end{equation}

We note that (with $c_1,c_2,\ldots$ as in \eqref{cn})
\[
\begin{split}
  \sqrt{I_2-\begin{pmatrix}\overline{a}\\ a\end{pmatrix}\begin{pmatrix}a&\overline{a}\end{pmatrix}}&=I_2-\sum_{n=1}^\infty c_n
  \left(\begin{pmatrix}\overline{a}\\ a\end{pmatrix}
    \begin{pmatrix}a&\overline{a}\end{pmatrix}
  \right)^n\\
  &=I_2-\sum_{n=1}^\infty c_n \left(\begin{pmatrix}\overline{a}\\ a\end{pmatrix}
    \begin{pmatrix}a&\overline{a}\end{pmatrix}\right)\left(\begin{pmatrix}\overline{a}& a\end{pmatrix}    \begin{pmatrix}a\\ \overline{a}\end{pmatrix}
  \right)^{n-1}\\
  &=I_2-\frac{\begin{pmatrix}\overline{a}& a\end{pmatrix}    \begin{pmatrix}a\\ \overline{a}\end{pmatrix}}{2|a|^2}\sum_{n=1}^\infty c_n(2|a|^2)^n\\
  &I_2-\frac{\begin{pmatrix}\overline{a}& a\end{pmatrix}    \begin{pmatrix}a\\ \overline{a}\end{pmatrix}}{2|a|^2}\left(1-\sqrt{1-2|a|^2}\right).
\end{split}
\]

\begin{theorem}
The function  $f$ belongs to $\mathfrak H(k)$ and  $f(a)=0$ if and only if
\[
  f(z)=b_a(z)g(z),
\]
with $g\in{\mathfrak  H(k)}^2$.
\end{theorem}

\begin{proof} We follow \cite{AT}. One direction is trivial while the converse is a direct consequence of \eqref{wer} with $g(z)=\begin{pmatrix}A_a^*f\\ B_a^*f\end{pmatrix}$
since $f(a)=0$.
\end{proof}

More generally, as in Proposition 4.5 and Section 5 of \cite{AT} we have
\begin{theorem}
  Let $z_1,\ldots, z_N\in B(0,1/\sqrt{2})$ and $w_1,\ldots, w_N\in\mathbb C$. There exists a Schur multiplier $s$ such that
  \begin{equation}
    s(z_n)=w_n,\quad n=1,\ldots, N
  \end{equation}
  if and only if the $N\times N$ matrix with $(n,m)$ entry
  \begin{equation}
\frac{1-w_n\overline{w_m}}{1-2{\rm Re}\, z_n\overline{z_m}}
   \end{equation}
  is non negative.
\end{theorem}
\begin{proof}
This holds thanks to the fact that the kernel $k(z,w)=\frac{1}{1-2{\rm Re}\, z\overline{w}}$ is a complete Nevanlinna-Pick kernel, and so the Nevanlinna-Pick interpolation problem is solved.
\end{proof}




\bibliographystyle{plain}

\end{document}